\documentclass[11pt]{article}

\usepackage{amsmath,amssymb,amsthm}
\usepackage{graphicx}
\usepackage{booktabs}
\usepackage{enumitem}
\usepackage{algorithm}
\usepackage{natbib}
\usepackage{array}    % for m{} column type
\usepackage{booktabs} % for \toprule, \midrule, \bottomrule, \addlinespace
\usepackage{pgf}
\usepackage{caption}
\usepackage{subcaption}
\usepackage[left=1in,right=1in,top=1in,bottom=1in]{geometry}
\usepackage{times}
\usepackage{tikz}
\usetikzlibrary{arrows.meta, positioning, calc, swigs, shapes}
\usetikzlibrary{arrows,shapes.arrows,shapes.geometric, shapes.multipartb,backgrounds,decorations.pathmorphing,positioning,fit,automata}
\tikzset{
  >=stealth',
  true/.style={
    rectangle, draw=black, very thick, text width=6.5em, minimum height=2em,
    text centered, fill=gray, opacity = 0.5},
  punkt/.style={
    rectangle, rounded corners, draw=black, very thick, text width=6.5em,
    minimum height=2em, text centered},
  est/.style={circle, draw=black, very thick, text centered},
  shade/.style={circle, draw=black, very thick, fill=gray!50, text centered},
  weight/.style={
    circle, draw=black, very thick, text width=6.5em, minimum height=2em, text centered},
  pil/.style={->, thick, shorten <=2pt, shorten >=2pt,},
  double/.style={<->, thick, shorten <=2pt, shorten >=2pt,},
  dash/.style={dashed, thick, shorten <=2pt, shorten >=2pt,},
  dashdouble/.style={<->, dashed, thick, shorten <=2pt, shorten >=2pt,}
}

\newcommand{\ind}{\perp\!\!\!\perp}

\newtheoremstyle{note}% <name>
{8pt}{8pt}{}{}{\bfseries}{:}{.5em}{}
\theoremstyle{note}

\newtheorem{remark}{Remark}

\newtheorem{definition}{Definition}

\newtheorem{proposition}{Proposition}
\allowdisplaybreaks

\begin{document}
\title{Causal Inference: A Tale of Three Frameworks}
   \date{}
    \maketitle

\begin{center}
    \vspace{-30pt}
%   \author{%\large Dingke Tang, Dehan Kong, and Linbo Wang \\
% %  \vspace{10pt}
%   Department of Statistical Sciences, University of Toronto}\\ \vspace{20pt}
%     \author{\large Dehan Kong \\
%     \vspace{10pt}
%   Department of Statistical Sciences, University of Toronto}\\ \vspace{10pt}

  \author{\large Linbo Wang$^{1}$, Thomas Richardson$^{2}$, and James Robins$^{3}$ \\
  \vspace{10pt}
$^{1}$Department of Statistical Sciences, University of Toronto \\
$^{2}$Department of Statistics, University of Washington \\
$^{3}$Harvard T.H. Chan School of Public Health}\\
  \vspace{30pt}
 {\large 
 }
\end{center}

% Abstract / Keywords per JDS template
\begin{abstract}
Causal inference is a central goal across many scientific disciplines. Over the past several decades, three major frameworks have emerged to formalize causal questions and guide their analysis: the potential outcomes framework, structural equation models, and directed acyclic graphs. Although these frameworks differ in language, assumptions, and philosophical orientation, they often lead to compatible or complementary insights. This paper provides a comparative introduction to the three frameworks, clarifying their connections, highlighting their distinct strengths and limitations, and illustrating how they can be used together in practice.  The discussion is aimed at researchers and graduate students with some background in statistics or causal inference who are seeking a conceptual foundation for applying causal methods across a range of substantive domains.
\end{abstract}


% ===== Main Content (your manuscript, unchanged except wrapper) =====

\section{Introduction}
Causal inference is the science of understanding the consequences of interventions, requiring assumptions that extend beyond those needed for purely associational analysis. Its importance has grown rapidly in the era of machine learning and artificial intelligence, where the ability to draw reliable causal conclusions is central to building systems that are not only predictive but also trustworthy, transparent, and robust to distributional shifts \citep{peters2016causal,Wachter2017HJLTCounterfactual,Pearl2019CACMSevenTools,Arjovsky2020IRM,Buhlmann2020StatSci,TjoaGuan2021TNNLSXAIHealth,Scholkopf2022CausalityML,Jiao2024ResearchSurvey}. Over the past decades, three foundational frameworks have emerged to formalize causal reasoning: the potential outcomes framework, nonparametric structural equation models (NPSEMs), and directed acyclic graphs (DAGs). Each framework carries its own formal machinery, conceptual underpinnings, and historical roots. Although they originated in distinct disciplinary traditions, they are now increasingly recognized as complementary, and in many cases translatable into one another.

A substantial literature surveys causal inference from within a single framework or with an emphasis on identification and estimation approaches \citep[e.g.][]{ImbensWooldridge2009JEL,spirtes2010introduction,pearl2010IJB,kuang2020causal,Yao2021Survey,li2023bayesian,Jiao2024ResearchSurvey}. However, there has been little concise, side-by-side treatment that translates assumptions and results across the three frameworks and clarifies when they agree and when they differ, despite earlier comparative discussions such as \cite{greenland2002overview}, \citet{pearl2009causality},  and \citet{robins2011alternative}. This review aims to fill that gap. We offer a selective review of the three frameworks, with the goal of elucidating their differences, similarities, and interconnections. Rather than offering a comprehensive survey, we aim to guide readers, particularly those new to the field or  approaching causal inference from a single perspective, through the conceptual landscape of causal inference using simple, illustrative examples.

In this review, we focus on \emph{interventional} questions, which ask what would happen under specific actions or treatments applied to the entire population, without conditioning on what actually occurred. These contrast with \emph{counterfactual} questions, which ask what would have happened for the same individual, had things been different,  given what actually occurred. For example, as noted by \citet[][Remark 6.20(ii)]{peters2017elements}, suppose someone offers you \$10,000 if you correctly predict the result of a coin flip. You guess ``heads'' and lose. The interventional effect of guessing ``heads'' versus ``tails'' is zero, whereas the counterfactual effect, given that you actually guessed ``heads'' and lost, is \$10,000. While counterfactuals play a central role in philosophical and legal reasoning, their identification typically requires stronger, often untestable, assumptions. See Remark~\ref{remark:counterfactual}  and  \citet{wu2026position} for further discussion of this distinction.

Although we do not delve into estimation methods or unmeasured confounding in depth, we briefly comment on these in later sections to clarify their connection to the foundational frameworks. Other approaches, such as Dawid’s decision-theoretic framework \citep{dawid2000causal,dawid2015statistical,richardson2023potential}, are also omitted from our discussion. The decision-theoretic perspective defines causality through conditional independence relations among decisions, observations, and consequences, without invoking counterfactuals or structural equations. Conceptually, it shares common ground with all three frameworks reviewed here: like potential outcomes, it distinguishes between observed and hypothetical actions; like structural models, it emphasizes invariance of mechanisms under intervention; and like causal DAGs, it encodes causal assumptions through independence relations. However, it differs from these frameworks in that it formulates causal reasoning entirely within the observable world, avoiding reference to unobservable counterfactual quantities. Because our focus is on frameworks that explicitly model potential outcomes or interventions, we do not explore the decision-theoretic approach in detail. For these topics, and for comprehensive treatments of any individual framework, we refer readers to existing textbooks \citep[e.g.,][]{pearl2009causality,peters2017elements,imbens2015causal,hernan2020causal} and the broader specialized literature.

Another influential framework that we do not review in detail is the sufficient–component cause (SCC) framework; see \cite{rothman1976causes}, \citet{greenland1988invariants}, \citet{greenland2002overview}, \citet[][Chap. 2]{rothman2008modern}, \citet{vanderweele2009minimal},
\citet[][\S7]{vanderweele2012general},
and \citet[][p.10 \& Chap. 10]{vanderweele2015explanation} for detailed discussions. SCC models provide a mechanistic representation of causation and are well suited for questions about explanatory causal structures and mechanistic interaction. These mechanistic questions differ somewhat from the interventional questions that motivate our review, and several SCC quantities cannot be expressed within the causal DAG formalism. Since our goal here is to compare three frameworks for interventional causal questions, we do not cover SCC models, although they remain an important component of the broader causal inference literature.

A further set of traditions, especially influential in health and legal applications, involves more informal approaches to causal inference alongside formal models. \cite{Greenland2004Overview}, for example, distinguishes canonical approaches based on considerations such as the Hill's list \citep{Hill1965Environment}, in which causality is treated as a property of an observed association to be judged using qualitative “symptoms” like temporality, strength, consistency, and dose response, from methodologic modeling or bias analysis, which uses parametric models and sensitivity analyses to assess the impact of selection, unmeasured confounding, and measurement error on effect estimates. These approaches are widespread in applied practice and can be combined with potential outcome, structural, and graphical models, but they fall outside our focus on three foundational frameworks for representing and identifying interventional effects.

The frameworks reviewed here primarily concern single-stage interventions. Many applied settings, however, involve sequential or adaptive decisions, such as treatment policies in medicine, dynamic pricing, or reinforcement learning. These settings have been extensively studied in the literature on dynamic treatment regimes and longitudinal causal inference, most notably through Robins’s development of g-methods and related approaches \citep{robins1986new,hernan2020causal}. These frameworks generalize the ideas discussed here to multi-stage interventions and time-varying treatments. Because our focus is on the conceptual connections among the three foundational frameworks, we do not cover these extensions in detail.

Our review proceeds as follows. We begin by introducing the three frameworks in turn, highlighting how each connects observational data to hypothetical interventions and how each encodes key assumptions, such as the absence of unmeasured confounding. In Section~\ref{sec:translating}, we explore how these frameworks can be formally related. We show how potential outcomes arise naturally within NPSEMs, and conversely, how systems of potential outcomes can be used to construct NPSEMs through canonical representations. We also examine how causal DAG models are implied by structural equation models under the assumptions of independent errors and autonomy, and how NPSEMs can, in turn, be generated from a given causal DAG using the functional representation lemma. Finally, we introduce Single World Intervention Graphs (SWIGs), which integrate the graphical and potential outcomes approaches by explicitly representing potential outcomes on a modified DAG, thereby bridging these frameworks without invoking 
all of the assumptions imposed by
%the full strength of 
structural models.

Section~\ref{sec:comparison} offers a comparative analysis of the three frameworks along several dimensions, including their philosophical orientation, expressive capacity, and identification power. We discuss when the additional assumptions imposed by NPSEMs with independent errors can lead to stronger identification results and when such assumptions may be overly restrictive. We conclude with practical guidance for applied researchers and reflect on recent methodological developments that seek to unify these perspectives into a coherent and transparent workflow for causal inference.

% We emphasize that this is a selective review. In particular, we note several important topics that are beyond our scope:
% (1) Beyond the three popular frameworks discussed here, other approaches to causal inference exist, such as the decision-theoretic framework developed by \citet{dawid2000causal,dawid2015statistical}.
% (2) Our focus is primarily on interventional questions. Counterfactual questions \citep[][p.~29]{pearl2009causality} differ in nature.   (3) We limit our discussion to introducing the three frameworks and do not delve into the rich literature on the identification and estimation of causal effects.

\section{Three Frameworks for Causal Inference}
\label{sec:three-frameworks}

The formalization of causal effects in mathematical terms is widely recognized as one of the major advances in the modern theory of causality. Unlike purely associational quantities, causal effects characterize how an outcome responds to different hypothetical interventions that set a treatment or exposure to alternative values. At a broader level, counterfactual reasoning considers what {\it would} have happened for the same individual under alternative scenarios that may be contrary to fact. This conceptual leap distinguishes causal inference from classical statistics, which has traditionally focused on describing correlations and conditional dependencies among \emph{observed} data, or on making inferences regarding the distribution from which a sample was drawn.

Causal inference can be understood through at least three principal frameworks: (1) Potential Outcomes, (2) SEMs, and (3) DAGs. While each framework stems from a different disciplinary tradition, all aim to bridge the gap between observational data and interventional scenarios. This section provides a detailed account of each framework, with the key assumptions summarized in Table \ref{tab:causal_frameworks}. We discuss their relationships and how they both unify and diverge in their approaches to causal reasoning in the following sections.

\begin{table}[ht]
    \centering
      \caption{How the three causal inference frameworks relate the observational and interventional worlds, and how they encode assumptions about unmeasured confounding. 
      Here, $A$ denotes the treatment, $Y$ denotes the outcome, and $L$ denotes  covariates }
\begin{tabular}{m{2cm} m{6cm} m{6cm}}
\toprule
Framework & How the observational world relates to the (hypothetical) interventional world & How to encode or check no unmeasured confounding \\
\midrule
Potential Outcomes 
  & Consistency: $Y = Y(A)$ 
  & Ignorability\textsuperscript{a}: $A \ind Y(a) \mid L$ \\
\addlinespace[0.6em]
NPSEM\textsuperscript{b} 
  & Autonomy\textsuperscript{c}: intervention modeled by replacing the structural equation for $A$, while keeping all other equations unchanged 
  & Independent errors\textsuperscript{d} (NPSEM-IE) \\
\addlinespace[0.6em]
Causal DAG 
  & Modularity: intervention modeled by removing incoming edges 
  & Causal sufficiency: no unmeasured common causes of variables in the graph; under this assumption, $L$ can be chosen to satisfy the backdoor criterion\textsuperscript{e} \\
\bottomrule
\end{tabular}
    \label{tab:causal_frameworks}
    \vspace{0.5em}

 \begin{flushleft}
    \footnotesize
       \textsuperscript{a} $\ind$ denotes independence.  \\[4pt]
        \textsuperscript{b} An NPSEM  typically assumes that interventions on all variables are well-defined.
        \\[4pt]
    \textsuperscript{c} In the NPSEM literature, this property is also called modularity, highlighting that the mechanism for $A$ can be replaced without altering other equations. \\[4pt]
    \textsuperscript{d} A generic NPSEM does not impose independent errors; the NPSEM with independent errors (NPSEM-IE) is stronger in that it adds this assumption.  \\[4pt]
    \textsuperscript{e} The scope differs across frameworks. In the potential outcomes framework, ignorability is typically imposed for a specific effect of $A$ on $Y$, and thus enables identification of that individual effect alone. In contrast, an NPSEM-IE encodes a full system of structural equations and independent errors, which enables identification of all the interventional effects in the model. Causal DAG models with causal sufficiency rule out hidden common causes; under this assumption, graphical tools such as the backdoor criterion provide effect-specific tests for identification across multiple treatment–outcome pairs.
    \end{flushleft}

\end{table}

\subsection{The Potential Outcomes Framework}

The potential outcomes framework provides an intuitive and mathematically tractable foundation for causal inference. It was originally introduced by \citet{neyman1923application} for randomized experiments and was later formalized and extended to observational studies by \citet{rubin1974estimating}. This was further extended by \cite{robins1986new} to studies with time-varying treatments and confounders, as well as to direct and indirect effects and feedback of one cause on another. This framework is closely related to an independent line of work in econometrics that developed the so-called switching regression model \citep{roy1951some, manski1999identification}.
 These developments allow researchers to pose well-defined causal questions, even in the absence of randomized experiments, by explicitly invoking assumptions that link observed data to potential outcomes. We refer interested readers to \citet{rubin2005causal, morgan2014counterfactuals,imbens2015causal,hernan2020causal} for a more detailed introduction.

Suppose the observed data consist of a triple \( (L, A, Y) \), where \( L \) denotes  (a vector of) covariates, \( A \in \{0,1\} \) is a binary treatment indicator (for example, placebo versus active treatment), and \( Y \) is the outcome of interest. The variables \( L \) capture characteristics   that may confound the treatment and outcome relationship. These covariates play a key role in adjustment procedures that aim to mimic randomization in observational settings.

For example, consider a job-training program $(A=1)$ intended to improve employment outcomes $(Y)$. Individuals with higher prior education or stronger motivation $(L)$ may be both more likely to enroll in the program and more likely to achieve better employment outcomes regardless of enrollment. This example provides a simple setting for comparing how the three causal frameworks represent and identify the same underlying causal effect.
 
For a binary treatment variable \( A \), we define two potential outcome variables:
\begin{itemize}
  \item \( Y(0) \): the value of \( Y \) that would be observed for a given unit if assigned \( A = 0 \) (placebo);
  \item \( Y(1) \): the value of \( Y \) that would be observed for a given unit if assigned \( A = 1 \) (treatment).
\end{itemize}
The variables \( Y(0) \) and \( Y(1) \) are two distinct random variables representing hypothetical outcomes under different treatment assignments. They are \emph{not} different realizations of the same variable.  For each unit, only one of these two outcomes is realized and observed, depending on the actual treatment received. The other remains counterfactual. This leads to the \textit{Fundamental Problem of Causal Inference} \citep{holland1986statistics}: we can never simultaneously observe both potential outcomes \( Y(0) \) and \( Y(1) \) for the same individual. From this perspective, the potential outcomes framework views causal inference as a missing data problem, since at least one of the two potential outcomes is missing for each unit.

To address this challenge, causal inference often focuses on population-level quantities that average over the distribution of potential outcomes. This leads to a distinction between two types of causal effects:\\[-8pt]

\begin{itemize}
\item \textbf{Individual causal effect}: Defined as \( Y_i(1) - Y_i(0) \), this measures the causal effect for a specific unit \( i \). Although conceptually straightforward, it is unobservable and generally not identifiable from data without strong assumptions or additional structure (such as deterministic effects or monotonicity).\\[-8pt]

\item \textbf{Average causal effect (ACE)}: Defined as
\(
\text{ACE} = \mathbb{E}\{Y(1) - Y(0)\},
\)
the ACE captures the expected difference in outcomes had everyone in the population received treatment versus control. This is the primary estimand (i.e. parameter of interest) in many empirical studies, and it is identifiable under assumptions such as complete randomization or conditional ignorability (see eqn. \eqref{eqn:weak-ignorability}).
\end{itemize}

\begin{remark}
   In addition to the individual and average causal effects, other causal estimands that offer more nuanced insights into treatment effects may also be of interest. One such quantity is the \textit{Conditional Average Treatment Effect (CATE)}, defined as \( \mathbb{E}\{Y(1) - Y(0) \mid X = x\} \), which represents the average treatment effect for individuals with covariates \( X = x \) where $X\subset L$. CATE is central to personalized decision-making and helps characterize treatment effect heterogeneity across subpopulations. Another commonly used estimand is the \textit{Effect of Treatment on the Treated (ETT)}, given by \( \mathbb{E}\{Y(1) - Y(0) \mid A = 1\} \), which captures the causal effect among those who actually received the treatment; the ETT is useful, for example, when evaluating whether a program is effective for those who are participating in it. The \textit{Quantile Treatment Effect (QTE)}, expressed as \( Q_{Y(1)}(\tau) - Q_{Y(0)}(\tau) \), quantifies the difference in the \( \tau \)-th quantiles of the potential outcome distributions under treatment and control. This allows researchers to assess how treatment affects the entire outcome distribution rather than just its mean. Importantly, in contrast to ACE, QTEs are not “individual-level” effects: the difference between the treatment and control quantiles need not correspond to the effect experienced by any particular individual. Instead, QTEs characterize population-level distributional shifts, which may provide valuable insight into distributional heterogeneity even if they do not map to individual units.
\end{remark}

\begin{remark}
  Potential outcomes are sometimes referred to as \emph{counterfactuals} in the literature. {\it Before} treatment assignment, both \(Y(1)\) and \(Y(0)\) are potential outcomes. {\it After} treatment assignment, the outcome corresponding to the received treatment is \emph{factual}, while the other is \emph{counterfactual}. (However, many papers use the terms potential outcomes and counterfactuals interchangeably.) To avoid ambiguity, we use the term \emph{potential outcomes} to collectively refer to both \(Y(1)\) and \(Y(0)\).
\end{remark}

\begin{remark}
\label{remark:counterfactual}
    The estimands we have introduced so far are sometimes referred to as \emph{interventional} estimands, as they pertain to the distribution of outcomes under hypothetical interventions, such as setting the treatment to a fixed value. Another type of question that may be of interest involves \emph{counterfactual} reasoning. For example, when \( Y \) is binary, the persuasion rate is defined as \( \mathbb{P}\{Y(1) = 1 \mid Y(0) = 0\} \) \citep{jun2023identifying}, which quantifies the probability that an individual would take an action under treatment given that they would not have done so under control. Other examples include the probability of necessity, \( \mathbb{P}\{Y(0) = 0 \mid A = 1, Y = 1\} \); the probability of sufficiency, \( \mathbb{P}\{Y(1) = 1 \mid A = 0, Y = 0\} \); and the probability of necessity and sufficiency, \( \mathbb{P}\{Y(1) = 1, Y(0) = 0\} \) \citep[\S~9.2]{pearl2009causality}.
 Identification of counterfactual estimands is typically much more involved, as it relies on the dependence between potential outcomes, which cannot be recovered from observed data, even in randomized controlled experiments. We refer readers to \citet[][Example 3.4]{peters2017elements} for an interesting example and discussions on whether counterfactuals or interventional quantities should serve as the basis for legal decisions; see also \citet{robins1989probability}, \citet{robins1991estimability},  \citet[][\S 1.4.4]{pearl2009causality}, and \citet{wu2026position}. Along with association questions, interventional and counterfactual questions are  referred to as the three levels of the ladder of causation by \citet[][\S~1]{pearl2018bookofwhy}.
\end{remark}

A crucial component of the  Potential Outcome  framework is the \textit{Stable Unit Treatment Value Assumption (SUTVA)} \citep{rubin1980randomization}, which posits that (i) there are no multiple versions of the treatment, and (ii) one unit's outcome is unaffected by the treatment assigned to another unit, meaning that there is no interference between units.

The first part of SUTVA is closely related to the \emph{consistency} assumption \citep{cole2009consistency,vanderweele2009concerning,pearl2010consistency}, which is the fundamental assumption that links the observational world to the potential outcomes in a hypothetical interventional world. The consistency assumption states that if an individual receives treatment level \( A = a \), then their observed outcome \( Y \) equals their potential outcome under that treatment:
\[
Y = Y(a) \quad \text{if } A = a, \quad \text{or equivalently,} \quad Y = Y(A).
\]
This identity presumes that the treatment is well-defined and administered uniformly across units. In other words, there are no multiple or ambiguous versions of treatment corresponding to the same value \( a \).  

 A classic example of an ill-defined treatment arises when the so-called ``treatment'' is the causal effect of obesity (defined as BMI $>30$) on heart  disease \citep{hernan2020causal}.
 
 %see for example the study \citep{millard2019casual}, 
 %which is determined by prior exposures like diet and physical activity and thus represents an intermediate outcome rather than a well-defined intervention . Suppose we define the ``treatment'' as having a BMI of 25. 
 % since, 
 % at present, there is no known agent that can causally affect BMI without that agent also causing the outcome through pathways that are not mediated by obesity. 
 %This value may correspond to very different behavioral or physiological profiles across individuals. One person may reach this BMI through diet alone, another through intensive exercise, and yet another through illness. These different pathways can lead to different outcomes even if the measured value of BMI is the same. Consequently, the treatment level does not correspond to a unique intervention, which violates the condition of no multiple versions and undermines the consistency assumption.

\begin{remark}
The consistency assumption discussed here is not the same as the notion of consistency used to describe a statistical property of an estimator \citep[e.g.,][]{lehmann2006theory}, which refers to convergence in probability to the true parameter value.
\end{remark}
 
The second part of SUTVA, namely the no interference assumption, implies that each unit's potential outcome depends only on its own treatment assignment and not on the treatment received by others. This assumption can be violated in many real-world contexts. For example, in vaccine studies, the likelihood of an individual contracting an infectious disease may depend not only on whether they are vaccinated but also on the vaccination status of others in their community. Such \emph{interference} is common in infectious disease epidemiology and social network settings, where spillover, peer effects, or herd immunity can induce dependencies across units.
A growing body of work has extended the potential outcomes framework to accommodate interference. We refer interested readers, for example, to \citet[][\S 9]{kuang2020causal} and \cite{bajari2023experimental} for surveys of this area.

To identify the ACE under the SUTVA, a key condition is the assumption of \emph{no unmeasured confounding}, also known as \emph{conditional (weak) ignorability} or \emph{selection on observables}:
\begin{equation}
\label{eqn:weak-ignorability}
    A \ind Y(a) \mid L, \quad a=0,1.
\end{equation}
 This condition implies that, conditional on \( L \), the treatment assignment \( A \) is as good as random and thus independent of the potential outcomes. Under this assumption and a positivity condition that ensures all levels of treatment are possible across values of \( L \), the ACE is identified by the following  formula:
\begin{equation}
\label{eqn:ace-identification}
\text{ACE} = \mathbb{E}_L \left\{ \mathbb{E}(Y \mid A = 1, L) - \mathbb{E}(Y \mid A = 0, L) \right\},
\end{equation}
where the outer expectation is over the marginal distribution of \( L \). Equation \eqref{eqn:ace-identification} was termed ``direct standardization'' in the early epidemiological literature \citep{neison1844method} and is also a special case of the g-formula of \citet{robins1986new}.
This expression can be estimated using regression, stratification, matching, inverse probability weighting, or doubly robust methods, depending on the structure and quality of the observed data.

\begin{remark}
    In contrast to \eqref{eqn:ace-identification}, the association between $A$ and $Y$ can be written as $E(Y\mid A=1) - E(Y\mid A=0) = \mathbb{E}_{L\mid A=1}  \mathbb{E}(Y \mid A = 1, L) -  \mathbb{E}_{L\mid A=0} \mathbb{E}(Y \mid A = 0, L).$
\end{remark}

\subsection{Structural Equation Models}

% A \emph{Nonparametric Structural Equation Model (NPSEM)}
An NSPEM \citep{haavelmo1943probability,strotz1960recursive,pearl2009causality,spirtes2000causation,halpern2000axiomatizing} consists of a finite set of random variables \( V = \{V_1, \ldots, V_p\} \), a collection of exogenous (unobserved) random variables \( \varepsilon = (\varepsilon_1, \ldots, \varepsilon_p) \), and a collection of measurable functions \( f_1, \ldots, f_p \), such that each endogenous variable \( V_j \) is determined by an equation of the form
\begin{equation}
\label{eqn:npsem-ie}
V_j = f_j(U_j, \varepsilon_j),     
\end{equation}
where \( U_j \subseteq V \setminus \{V_j\} \) denotes a subset of variables that serve as inputs to the structural assignment of \( V_j \), known as the parents of $V_j$. In this context, the variables in \( V \) are referred to as \emph{endogenous} variables because their values are determined within the system by the structural equations. In contrast, the components of \( \varepsilon \) are called \emph{exogenous} variables, as they represent factors external to the system whose values are not influenced by other variables in the model. These exogenous variables typically encode latent background factors or sources of randomness that affect the system but are not explained by it.

Together, the functions \( f_j \) and exogenous variables \( \varepsilon_j \) specify a data-generating process 
and thereby induce a joint distribution over \( V \). 
Importantly, NPSEMs do not need to be acyclic: cycles or feedback relations can in principle be allowed. 
When the system is \emph{acyclic}, however, there exists (at least) one ordering 
\( V_{\pi(1)}, \ldots, V_{\pi(p)} \) such that each \( V_{\pi(j)} \) depends only on variables earlier in the ordering. 
In this case, the structure of dependence among variables can be represented by a directed acyclic graph (DAG) 
(see Section~\ref{sec:npsemtodag}). 
The model is then fully defined by the tuple \( (f_1, \ldots, f_p, \varepsilon_1, \ldots, \varepsilon_p) \) 
together with a valid causal ordering (e.g., as obtained from Simon's causal ordering algorithm \citep{simon1953causal}). 
The nonparametric nature of NPSEMs arises from the fact that no parametric form is assumed for the functions \( f_j \), 
allowing for considerable flexibility in modeling complex causal relationships. 
A key assumption in many NPSEM applications is that the exogenous variables in \( \varepsilon \) are jointly independent, 
yielding the special case known as the NPSEM with Independent Errors (NPSEM-IE).  

As an example, in a causal system involving a covariate \( L \), a treatment \( A \), and an outcome \( Y \), an NPSEM may specify the following equations:
\begin{align}
L &= f_L(\varepsilon_L), \notag \\
A &= f_A(L, \varepsilon_A), \label{eq:npsem} \\
Y &= f_Y(L, A, \varepsilon_Y), \notag 
\end{align}
where $ \varepsilon_L, \varepsilon_A,$ and $\varepsilon_Y $ are exogenous error terms representing unmeasured factors.%
In the job-training example, these equations can be interpreted as representing how individual characteristics $(L)$ influence program participation $(A)$ and, in turn, employment outcomes $(Y)$. This formulation expresses the data-generating process in a way that makes causal relations explicit, distinguishing between background factors and the mechanisms that connect them.

\begin{remark}
    The {\em Causally Interpretable Structured Tree Graph (CISTG) as detailed as the data} \protect\citep{robins1986new} is mathematically equivalent to an NPSEM model (with no assumption on the errors); see  Equations (\ref{eq:rewrite}) and (\ref{eqn:npsem-canonical}) below.
\end{remark}

On the surface, the NPSEM in \eqref{eq:npsem} appears identical to a nonparametric regression model for \( (L, A, Y) \). 
However, unlike statistical regression models that describe relationships in the observational world only, 
what makes NPSEMs \emph{structural} is that they also define what would happen in an interventional world by 
removing or replacing any subset of structural equations. 

The key assumption in the SEM framework that links the observational world to interventional worlds 
is the \emph{autonomy assumption} \citep{aldrich1989autonomy}, also known as the \emph{invariance} or 
\emph{modularity} assumption \citep{peters2017elements}. 
It asserts that each structural equation represents an independent, self-contained causal mechanism 
that remains invariant under interventions on other variables. 
Formally, an \emph{intervention} on a variable \(A\) corresponds to removing the structural equation for \(A\) and replacing it with a constant assignment. This is denoted by Pearl’s \emph{do-operator} \citep{pearl2009causality}, written as \(\text{do}(A=a)\).  Historically, \citet{robins1986new} used the notation \(g = a\) to express the same idea; Pearl's do notation is now more widely used.
%although the \(g\)-notation is no longer standard and has since been replaced by Pearl's \(\text{do}\)-operator.
% footnote changed

In the job-training example, this corresponds to evaluating a hypothetical scenario in which program participation is externally assigned rather than self-selected. The NPSEM in \eqref{eq:npsem} implies that in an interventional world where \( A \) is fixed to \( a \), 
the system becomes:
\begin{align*}
L &= f_L(\varepsilon_L), \\
A &= a, \\
Y &= f_Y(L, a, \varepsilon_Y).
\end{align*}

Under this formulation, the ACE can be represented as
%\footnote{This ACE can also be represented as $E[Y(a=1) - Y(a=0)]$
%because 
%under the NPSEM (\ref{eq:npsem}
%$f_Y(L, a, \varepsilon_Y) =  Y(a)$; see 
%  Equations (\ref{eq:rewrite}) and (\ref{eqn:npsem-canonical}) below.}
\[
\text{ACE} = \mathbb{E}\{Y \mid \text{do}(A = 1)\} - \mathbb{E}\{Y \mid \text{do}(A = 0)\} 
=  \mathbb{E}\{ f_Y(L,1,\varepsilon_Y) - f_Y(L,0,\varepsilon_Y) \}.
\]

Autonomy ensures that interventions can be represented by replacing the structural equation for the intervened variable, without having to re-specify the rest of the system.
This feature allows us to predict counterfactual outcomes under hypothetical interventions. 
In real-world systems, however, autonomy is not always guaranteed \citep[e.g.,][]{cartwright2007hunting}. 
Administering the treatment may change the patient's behavior or alter how the outcome is measured. 
For instance, it may modify the reporting mechanism or induce side effects 
that feed back into the outcome-generation process in an unmodeled way. 
In such cases, the function \( f_Y \) may no longer be invariant to interventions on \( A \). 
In other words, the mapping from \( (L, A, \varepsilon_Y) \) to \( Y \) may differ depending on 
whether \( A \) is set naturally or via intervention.

\begin{remark}
   The term Structural Causal Models (SCMs) is often used in the literature. It may refer to the NPSEM \citep[e.g.,][p.~203, Definition~7.1.1]{pearl2009causality} alone or specifically to the NPSEM with independent errors (NPSEM-IE) \citep[e.g.,][p.~44, Definition~2.2.2]{pearl2009causality}; see also \citet[][p.~62, Definition~6.2]{peters2017elements}. Following \citet{shpitser2022multivariate}, we use our terminology here to emphasize the distinction between NPSEM and NPSEM-IE.
\end{remark}

\begin{remark}
The independent error and autonomy assumptions in the NPSEM-IE framework 
imply the \emph{principle of independent mechanisms (PIM)} 
\citep{janzing2010causal,peters2017elements}. 
This principle states that the causal generative process can be decomposed into 
independent modules, one for each variable, that remain invariant to changes in the 
distributions of other variables. 
In the special case of two variables, PIM reduces to the so-called 
\emph{independence of cause and mechanism}, which asserts that the conditional distribution 
\(P(Y \mid X)\) is independent of the marginal distribution \(P(X)\). 
In other words, under NPSEM-IE assumptions, the causal mechanism mapping 
\(X \mapsto Y\) is invariant to changes in the distribution of \(X\).
\end{remark}

% \begin{remark}
%  \todo   \tblue{Add a remark on causal invariance that comes out of the autonomy assumption. Also on the result from anchor regression eqn. (1).} 
% \end{remark}

% \bigskip
% \bigskip
% \bigskip

% One conceptual distinction between NPSEMs and the potential outcomes framework is the direction of modeling: whereas potential outcomes begin with counterfactual definitions and derive observed data through consistency, NPSEMs begin with structural equations for observed variables and derive counterfactuals by removing and modifying equations. This difference can be subtle but has important implications for model specification and interpretation.

\subsection{Directed Acyclic Graphs}
\label{sec:dag}

The use of DAGs for causal inference stems from foundational work by \cite{spirtes2000causation} and \cite{pearl1995causal}, building on earlier ideas from Sewall Wright's path diagrams \citep{wright1921correlation}. In this framework, a causal system is represented by a DAG, where nodes correspond to variables and directed edges encode possible direct causal effects. In particular, the absence of an edge asserts that no direct causal effect is present. In the following, we introduce the DAG model from three perspectives: purely graphical, statistical/probabilistic, and causal.

\subsection*{Directed Acyclic Graph}

A \emph{Directed Acyclic Graph} (DAG) is a finite set of nodes connected by directed edges such that no cycles are present. 
In other words, there is no way to start at a node and return to it by following a sequence of directed edges.  

Formally, let \( \mathcal{G} = (V, E) \) be a DAG, where
\begin{itemize}
    \item \( V = \{V_1, \ldots, V_p\} \) is a set of vertices, and
    \item \( E \subseteq V \times V \) is a set of directed edges.
\end{itemize}

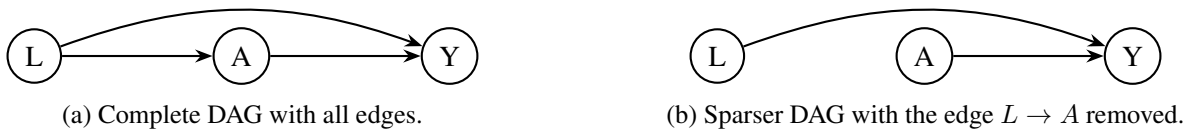
\begin{figure}[h]
    \centering

    % Subfigure (a): Complete DAG
    \begin{subfigure}[b]{0.45\textwidth}
        \centering
        \begin{tikzpicture}[node distance=2cm, >=Stealth, thick]
            \node[circle, draw] (L1) {L};
            \node[circle, draw, right=of L1] (A1) {A};
            \node[circle, draw, right=of A1] (Y1) {Y};

            \draw[->] (L1) -- (A1);
            \draw[->] (L1) to[out=20, in=160] (Y1);
            \draw[->] (A1) -- (Y1);
        \end{tikzpicture}
        \caption{Complete DAG with all edges.}
        \label{fig:dag_complete}
    \end{subfigure}
    \hfill
    % Subfigure (b): Given DAG
    \begin{subfigure}[b]{0.45\textwidth}
        \centering
        \begin{tikzpicture}[node distance=2cm, >=Stealth, thick]
            \node[circle, draw] (L2) {L};
            \node[circle, draw, right=of L2] (A2) {A};
            \node[circle, draw, right=of A2] (Y2) {Y};

            \draw[->] (A2) -- (Y2);
            \draw[->] (L2) to[out=20, in=160] (Y2);
            % No edge from X2 to Y2
        \end{tikzpicture}
        \caption{Sparser DAG with the edge \( L \rightarrow A \) removed.}
        \label{fig:dag_sparse}
    \end{subfigure}
    \caption{Two example DAGs over variables \( L \), \( A \), and \( Y \).}
    \label{fig:dag_example}
\end{figure}

Figure~\ref{fig:dag_example} shows two example DAGs. We next introduce some graph terminology that will be useful later. At this point, no statistical or causal interpretation is attached to the graph.

In a DAG $\mathcal{G}$, the \emph{parents} of a node \( V_j \), denoted $\text{Pa}_\mathcal{G}(V_j)$, or \( \text{Pa}(V_j) \) in short, are all nodes with arrows pointing directly into \( V_j \):
\[
\text{Pa}(V_j) = \{ V_k \in V : (V_k \to V_j) \in E \}.
\]
For example, in both Figure~\ref{fig:dag_complete} and Figure~\ref{fig:dag_sparse}, \( \text{Pa}(Y) = \{L, A\} \). If \( V_k \) is a parent of \( V_j \), then
$V_j$ is a \emph{child} of $V_k.$

A \emph{path} between nodes \(V_k\) and \(V_j\) is a sequence of distinct nodes
\[
V_k - V_{i_1} - V_{i_2} - \cdots - V_j
\]
such that each consecutive pair is connected by an edge, regardless of its orientation. 
A \emph{directed path} or \emph{causal path} from \(V_k\) to \(V_j\) is a path in which all edges 
point forward along the sequence, for example 
\(V_k \to V_{i_1} \to \cdots \to V_j\).  
Any other path from \(V_k\) to \(V_j\) is a \emph{non-causal path}. 
A particularly important class of non-causal paths are the 
\emph{backdoor paths}, which are paths from \(V_k\) to \(V_j\) that begin with an arrow into \(V_k\). 

The \emph{ancestors} of \( V_j \), denoted \( \text{An}(V_j) \), are all nodes from which \( V_j \) is reachable by a directed path, with \( V_j \) itself included as an ancestor  of $V_j$.  
Dually, the \emph{descendants} of \( V_j \), denoted \( \text{De}(V_j) \), are all nodes reachable from \( V_j \) by a directed path, again including \( V_j \) itself. 
The set of \emph{non-descendants} of \( V_j \) is
\[
\text{ND}(V_j) := V \setminus \text{De}(V_j).
\]
For example, in Figure~\ref{fig:dag_complete} we have, for node \( L \),
\(
\text{An}(L)=\{L\},\ \text{De}(L)=\{L,A,Y\},\ \text{ND}(L)=\varnothing.
\)
In the sparser DAG of Figure~\ref{fig:dag_sparse}, for the same node \( L \),
\(
\text{An}(L)=\{L\},\ \text{De}(L)=\{L,Y\},\ \text{ND}(L)=\{A\}.
\)

Lastly, let \( A, B, C \) be disjoint subsets of nodes in a DAG \( \mathcal{G} \).   A path between a node in \( A \) and a node in \( B \) is said to be \emph{blocked} by \( C \) if it contains a subpath of one of the following forms:
\begin{enumerate}
    \item A \emph{chain} \( X \to M \to Y \),
    \( X \leftarrow M \leftarrow Y \)
    or a \emph{fork} \( X \leftarrow M \to Y \) such that \( M \in C \);
    \item A \emph{collider} \( X \to M \leftarrow Y \) such that neither \( M \) nor any of its descendants is in \( C \).
\end{enumerate}
Some authors refer to the  subpaths in 1. as {\em non-colliders}.
If all paths from nodes in \( A \) to nodes in \( B \) are blocked by \( C \), then \( A \) and \( B \) are said to be \emph{d-separated} by \( C \) (the \emph{d} denotes \emph{directional} \citep[][p. 46]{pearl2016causal}). Equivalently, d-separation can be defined using the notion of ``moral graph'' \citep{lauritzen1990independence}.
% \footnote{%
% The term ``ancestral graph'' here follows \citet{lauritzen1990independence}, where it denotes the induced subgraph on a set of nodes and their ancestors. 
% In later literature, particularly \citet{richardson2002ancestral}, the same term is used for a more general graphical object that may include bi-directed edges. 
% We use the \citet{lauritzen1990independence} definition throughout this paper.}
Specifically, consider the subgraph induced by the ancestors of \( A \cup B \cup C \). For each pair of non-adjacent nodes that share a common child, add an undirected edge between them (thus ``marrying'' the parent variables). Then, remove the directions of all edges to obtain an undirected moral graph. We say \( A \) and \( B \) are d-separated by \( C \) in the original DAG if and only if every path between \( A \) and \( B \) is intercepted by a node in \( C \) in the resulting moral graph.

Unlike the concepts of parents and descendants, \emph{d-separation} is a relatively modern development \citep{pearl1986graphoids, pearl1988probabilistic, geiger1990identifying,pearl1995causal}. It enables one to read off conditional independencies from directed graphs in probabilistic DAG models, to which we now turn.

\subsection*{Probabilistic DAG Model}

In probabilistic or statistical modeling, a DAG is sometimes referred to as a Bayesian network \citep{pearl1985bayesian}, and each node in a DAG represents a random variable. Here the DAG is 
used to encode the conditional independence structure of a joint distribution over a set of random variables \( (V_1, \ldots, V_p) \). Specifically, let \( P \) be a joint distribution over \( (V_1, \ldots, V_p) \), and let \( \mathcal{G} \) be a DAG. We say that \( P \) is \emph{Markovian} with respect to \( \mathcal{G} \) if it satisfies the following properties.

\begin{definition}[Markov Properties]
\label{def:markov}
\quad 
    \begin{enumerate}
        \item \textbf{Markov Factorization Property:} The joint distribution factorizes as:
        \begin{equation}
            P(V_1, \ldots, V_p) = \prod_{j=1}^p P(V_j \mid \text{Pa}_\mathcal{G}(V_j)). \label{eqn:factorization}
        \end{equation}
         
        \item \textbf{Global Markov Property:} For any disjoint sets \( A, B, C \subseteq V \), if \( C \) d-separates \( A \) and \( B \) in \( \mathcal{G} \), then
        \(
        A \ind B \mid C   \text{ in } P.
        \)      

        \item \textbf{Local Markov Property:} Each variable is conditionally independent of its non-descendants given its parents:
        \(
        V_j \ind \text{ND}(V_j) \mid \text{Pa}(V_j).
        \)
    \end{enumerate}
\end{definition}
These Markov properties are equivalent if the joint distribution of \( V \) has a density with respect to a product measure  \citep[][Theorem 6.22] {lauritzen1996graphical}.

Due to the global Markov property, in probabilistic DAG models, d-separation serves as a bridge between the structure of a DAG and the independencies encoded in the corresponding distribution. If the DAG is incomplete (i.e., contains missing edges), Pearl's \emph{d-separation criterion} can be used to identify the conditional independence relationships implied by the factorization.

As an example, consider the DAG shown in Figure~\ref{fig:dag_sparse}. The DAG encodes the following factorization of the joint distribution as in \eqref{eqn:factorization}:
\(
P(L, A, Y) = P(L) \cdot P(A) \cdot P(Y \mid  L, A).
\)
 The local Markov property implies that \( A \) is independent of its non-descendant \( L \).
    On the other hand, using d-separation, we see that \( A \ind L \), because the only path connecting them, \( A  \rightarrow Y \leftarrow L \), contains a collider \( Y \).

\begin{remark}
\label{remark:markov}
Different graphs may encode the same set of conditional independence relationships. For example, the graphs \( A \rightarrow B \rightarrow C \) and \( A \leftarrow B \rightarrow C \) both encode the assumption that \( A \ind C \mid B \). In this case, these graphs are called \textit{Markov equivalent}. Graphs that are Markov equivalent cannot be distinguished based solely on conditional independencies in observational data. However, they may be distinguished in the presence of additional assumptions about the joint distribution, such as constraints on the functional form of the conditional distributions.  For detailed discussions, see  Example 3 in Section \ref{example:3} and \citet[\S 4 \& \S 7]{peters2017elements}. 
\end{remark}
We note that, in general, the converse of the (global) Markov property does not necessarily hold. 
That is, a probability distribution $P$ that is Markovian with respect to a DAG  may exhibit additional conditional independence relations 
that are not entailed by d-separation in the DAG. 
Such situations arise, for example, when special parameter values induce cancellations 
or deterministic relationships, leading to additional independencies not implied by the graph via d-separation. 
When the converse does hold, meaning that every conditional independence in $P$
corresponds exactly to a d-separation in $\mathcal{G}$, we say that the distribution is
\emph{faithful} to the graph.

\begin{definition}[Faithfulness]
\label{def:faithfulness}
Given a DAG $\mathcal{G}$ and a joint probability distribution $P$ 
over its variables, $P$ is \emph{faithful} to $\mathcal{G}$ if:
\begin{equation*}
A \ind_P B \mid C 
\quad \Leftrightarrow \quad 
A \text{ and } B \text{ are d-separated by } C \text{ in } \mathcal{G}.
\end{equation*}
\end{definition}

Not all probability distributions can be faithfully represented by a DAG. 
In particular, the graph structure enforces certain logical closure properties 
for d-separation that need not hold for conditional independence relations in general distributions. 
For example, if $C$ d-separates $A$ from $B$ and also $A$ from $D$, 
then $C$ must d-separate $A$ from $(B,D)$ in the same DAG. 
 However, this rule does not generally hold for probabilistic conditional independence. A distribution may satisfy both $A \ind_P B \mid C$ and $A \ind_P D \mid C$ without satisfying $A \ind_P (B,D) \mid C$.

It can be shown that the set of distributions unfaithful to a graph $\mathcal{G}$ has Lebesgue measure zero \citep[][Theorem 3.2]{spirtes2000causation}; see also \cite{uhler2013geometry}. This result, however, comes with an important caveat: there exist sequences of faithful distributions that converge arbitrarily close to unfaithful distributions. A more robust alternative is the concept of $\lambda$-strong faithfulness, which requires that the dependence strength (e.g., partial correlation) between d-connected variables exceeds some threshold $\lambda > 0$. Unlike exact faithfulness, the set of distributions that are not $\lambda$-strong faithful  may occupy a substantial volume in the space of distributions \citep{uhler2013geometry}.

A careful reader might observe that the first paragraph after Definition~\ref{def:faithfulness} seems to suggest that, when sampling uniformly over all distributions, the faithfulness assumption fails with probability one.
In contrast, the second paragraph implies that faithfulness holds almost surely. 
To reconcile this apparent contradiction, consider three binary variables $A$, $B$, and $D$, 
whose joint distribution is specified by probabilities $p(a,b,d) \in [0,1]$ for $(a,b,d) \in \{0,1\}^3$, 
subject to the constraint $\sum_{a,b,d} p(a,b,d)=1$. 
The space of all such distributions is the 7-simplex
\[
\Delta_7 \;:=\; \Bigl\{\, p \in [0,1]^8 : \sum_{a,b,d} p(a,b,d)=1 \Bigr\},
\]
which is a $7$-dimensional affine subset of $\mathbb{R}^8$.

Now consider the following subsets of distributions:
\begin{flalign*}
    \mathcal{S}_1 &= \{P_{A,B,D} : A \ind B\}; \\
    \mathcal{S}_2 &= \{P_{A,B,D} : A \ind B, \,\;\; A \ind D\}; \\
    \mathcal{S}_3 &= \{P_{A,B,D} : A \ind (B, D)\}.
\end{flalign*}
Here $\mathcal{S}_1$ corresponds to the distributions obeying the Markov property for the DAG ($B \rightarrow D \leftarrow A$),
while $\mathcal{S}_3$ gives the set of distributions Markov with respect to the DAG ($B\rightarrow D \quad A$). 
The sets $\mathcal{S}_1$, $\mathcal{S}_2$, $\mathcal{S}_3$ have dimensions 6, 5, and 4, respectively, 
within the 7-simplex $\Delta_7$.

Now suppose we  sample uniformly from the set $\mathcal{S}_1$. Then the subset of such distributions that are unfaithful to the DAG ($B\rightarrow D \leftarrow A$), 
such as those in $\mathcal{S}_2$, has measure zero. Thus, faithfulness holds with probability one 
under this sampling scheme.

In contrast, if we sample uniformly
 from $\mathcal{S}_2$, the set of distributions that satisfy 
$A \ind B$ and $A \ind D$, then with probability zero the resulting 
distribution will also satisfy $A \ind (B, D)$, i.e., fall into the set $\mathcal{S}_3$. 
Therefore, the distributions in $\mathcal{S}_2$ cannot be faithfully represented by any DAG almost surely. 
%In graphical terms, the model $\mathcal{S}_2$ cannot be represented by a DAG because,
%if $A$ is d-separated from $B$, and $A$ is d-separated from $D$, then $A$ is d-separated from $\{B,D\}$.

\begin{remark}
Strictly speaking, the phrase ``sample uniformly'' is not well defined without first specifying a parameterization, a point related to \emph{Borel's paradox}. Conditional distributions on measure zero sets can depend on the chosen parameterization \citep[][p.~441, Problem~33.1]{billingsley1995probability}; see also \cite{wang2022homogeneity}. In particular, subsets such as $\mathcal S_1$ and $\mathcal S_2$ are lower dimensional manifolds within the probability simplex, so there is no canonical uniform measure on them. For $\mathcal S_1$ and the associated DAG $B \rightarrow D \leftarrow A$, we make this precise by drawing each conditional probability table independently from a $\text{Dirichlet}(1,\ldots,1)$ prior and combining them using the DAG factorization, under which the set of unfaithful distributions has Lebesgue measure zero. For $\mathcal S_2$, we describe in Appendix~A1 a sampling scheme under which a distribution satisfying $A \perp B$ and $A \perp D$ almost surely does not satisfy $A \perp (B, D)$.
\end{remark}

% \tblue{Reconcile this difference.}

%% To resolve: Has measure 0 or 1?

\subsection*{Causal DAG Model}

So far, we have described how DAGs encode \emph{statistical} conditional independencies via the factorization of the observed data distribution. To endow a DAG with causal meaning, we must go beyond this statistical interpretation and commit to assumptions about how the system behaves under external interventions \citep[][Definition 1.3.1]{pearl2009causality}; see also \citet[][Definition 6.32]{peters2017elements}.

\begin{definition}[Causal DAG Model] \label{def:causal-dag}
A causal DAG model over random variables \(V=\{V_1,\ldots,V_p\}\) consists of a directed acyclic graph \(\mathcal{G}\) and a family of interventional distributions \(\{P_{\text{do}(v_S)}\}\) such that:
\begin{enumerate}
\item \textbf{(Observational Markov condition)} The \emph{observational distribution} \(P\) factorizes as \eqref{eqn:factorization}.

\item \textbf{(Modularity)} 
% For any atomic intervention \(do(V_j=v_j)\) (here \(V_j\) is the random variable and \(v_j\) the fixed value),
% \begin{equation}
% P_{do(V_j=v_j)}(V_1,\ldots,V_p)
% = \delta(v_j)\;\prod_{k\neq j} P\!\bigl(V_k \mid \mathrm{Pa}_{\mathcal{G}}(V_k)\bigr),
% \label{eqn:truncated-atomic}
% \end{equation}
% where \(\delta(v_j)\) denotes a point mass at \(v_j\). The \emph{mutilated} graph
% \(\mathcal{G}_{do(V_j=v_j)}\) is obtained by deleting all incoming edges into \(V_j\); 
% \emph{outgoing} edges from \(V_j\) are retained. Thus for \(k\neq j\), the parent set 
% \(\mathrm{Pa}_{\mathcal{G}}(V_k)\) is the same as in \(\mathcal{G}\), and if \(V_j\in \mathrm{Pa}_{\mathcal{G}}(V_k)\) 
% the factor \(P(V_k\mid \mathrm{Pa}_{\mathcal{G}}(V_k))\) is evaluated with \(V_j=v_j\).
For any joint intervention on a set \(S\subseteq [p] \equiv \{1,\ldots,p\}\) that sets $V_S$ to a fixed value $v_S$:
%\begin{equation}
%P_{\text{do}(V_S=v_S)}(V_1,\ldots,V_p)
%= \Bigl(\prod_{s\in S}I(V_s = v_s)\Bigr)\;
%  \prod_{k\notin S} P\!\bigl(V_k \mid \mathrm{Pa}_{\mathcal{G}}(V_k)\bigr),
%\label{eqn:truncated-distribution}
%\end{equation}
\begin{equation}
P_{\text{do}(v_S)}(V_{[p]\setminus S})
= %\Bigl(\prod_{s\in S}I(V_s = v_s)\Bigr)\;
  \prod_{k\notin S} P\!\bigl(V_k \mid V_{\mathrm{Pa}_{\mathcal{G}}(V_k)\setminus S}, V_{\mathrm{Pa}_{\mathcal{G}}(V_k)\cap S} = v_{\mathrm{Pa}_{\mathcal{G}}(V_k)\cap S} \bigr).
\label{eqn:truncated-distribution}
\end{equation}
% where the graph \(\mathcal{G}_{do(V_S=v_S)}\) is obtained by deleting all incoming edges into each \(V_s\) with \(s\in S\),
% while all outgoing edges from variables in \(S\) remain. Equivalently, for each \(k\notin S\) the factor 
% \(P(V_k\mid \mathrm{Pa}_{\mathcal{G}}(V_k))\) is evaluated with the parents in \(S\) clamped to their interventional values \(v_S\).
\end{enumerate}
\end{definition}

% \tblue{Add remarks about Markov equivalence (pearl p.~19) and causal discovery.}

To illustrate the modularity assumption in Definition~\ref{def:causal-dag}, consider Figure~\ref{fig:dag_complete}. When interpreted as a probabilistic DAG, this graph corresponds to a full model on the observed data distribution over \( (L, A, Y) \), since no edges are missing. However, if we interpret the same DAG causally, then the effect of an intervention on \(A\) 
is captured by the truncated DAG in Figure~\ref{fig:dag_sparse}, where all incoming edges into \(A\) 
have been removed. Formally,  \eqref{eqn:truncated-distribution} yields
\begin{equation}
    \label{eqn:truncated}
    P_{\text{do}(a)}(L,Y)
\;=\;
P(L)\, P(Y \mid A=a, L),
\end{equation}
which can be obtained by removing the term \( P(A \mid L) \) from the observational factorization 
and fixing \( A \) at the value \( a \). 
Marginalizing over   \(L\) gives
\begin{equation}
P_{\text{do}(a)}(Y) \;=\; \sum\limits_{L=l}P(L=l)\, P(Y \mid A=a,L=l),    
\end{equation}
which coincides with the g-formula   \eqref{eqn:ace-identification}.
This truncated factorization forms the basis of Pearl's \emph{do-calculus}, introduced shortly before Remark~\ref{remark:id}.

\begin{remark}
    Eqn. \eqref{eqn:truncated-distribution} was originally
    introduced by
    \citet{robins1986new}
and termed the g-computation algorithm formula  or simply the g-formula. It was subsequently independently re-discovered and named the manipulation theorem \citep{spirtes2000causation} and   the truncation formula \citep{pearl2009causality}.
\end{remark}

When a causal DAG is postulated, it is often assumed that
%An assumption commonly made in a causal DAG model is that 
the set of variables \( V \) is \emph{causally sufficient}, in the sense that there are no unmeasured common causes \( C \notin V \) that cause two or more variables in \( V \) \citep{spirtes2010introduction}. This assumption may be slightly relaxed to the notion of \emph{interventional sufficiency}, which is sufficient to ensure that the truncation formula~\eqref{eqn:truncated-distribution} holds; see \citet[\S~9.1]{peters2017elements} for a detailed discussion. When unmeasured common causes do exist between variables, it is conventional to 
represent them by adding bidirected edges (\(\leftrightarrow\)) between the affected 
nodes, resulting in an Acyclic Directed Mixed Graph (ADMG)~\citep{richardson2003markov}. 
Formally, such graphs can be obtained as \emph{latent projections} of DAGs with 
unobserved variables~\citep{verma1991equivalence, richardson2002ancestral}.

\subsection*{Definition and Identification of Causal Effects}

Under the causal DAG model, the \emph{average causal effect} (ACE) of a binary treatment \( A \in \{0, 1\} \) on an outcome \( Y \) is defined as:
\[
\text{ACE} = \mathbb{E}_{\text{do}(a=1)}(Y) - \mathbb{E}_{\text{do}(a=0)}(Y),
\]
where the expectation is taken with respect to the interventional distributions $P_{\text{do}(a)}(Y)$ defined via the truncated factorization in eqn. \eqref{eqn:truncated-distribution}.

In general, the interventional quantity \( \mathbb{E}_{\text{do}(a)}(Y) \) is not equal 
to its observational counterpart \( \mathbb{E}(Y \mid A = a) \) because of the presence of confounding. 
To recover causal effects, we may seek an \emph{adjustment set} \(L\) of observed variables 
such that conditioning on \(L\) blocks all non-causal paths from \(A\) to \(Y\).

When all variables in the DAG are observed, it is sufficient to include all the parents of \(A\) 
in the adjustment set \(L\), in which case the average causal effect (ACE) is identified by 
eqn.~\eqref{eqn:ace-identification} \citep[Theorem~3.2.2]{pearl2009causality}. 
More generally, even if some variables in the DAG are unobserved, the following 
\emph{backdoor criterion}  provides a graphical condition under which causal effects 
can still be identified from observational data.

A set of variables \( L \) satisfies the \emph{backdoor criterion} relative to treatment \( A \) and outcome \( Y \) if:
\begin{enumerate}
    \item No node in \( L \) is a descendant of \( A \), and
    \item \( L \) blocks all backdoor paths from \( A \) to \( Y \) in the DAG.
\end{enumerate}
Under this criterion, the ACE is still identified by eqn. \eqref{eqn:ace-identification} \citep[Theorem~3.3.2]{pearl2009causality}.

In the job-training example, background factors such as education and motivation (\(L\)) may influence both participation in the program (\(A\)) and employment outcomes (\(Y\)). In a corresponding DAG, these relationships are represented by arrows \(L \to A\) and \(L \to Y\), together with a causal arrow \(A \to Y\); see, for example, Figure~\ref{fig:dag_complete}. If \(L\) blocks all backdoor paths from \(A\) to \(Y\) and contains no descendants of \(A\), as is the case in Figure~\ref{fig:dag_complete}, then adjustment for \(L\) satisfies the backdoor criterion and identifies the causal effect of program participation.

\begin{remark}
   Although \citet[][\S~3.2.1, eqn.~(3.2)]{pearl2009causality} interprets the causal semantics of a DAG using the NPSEM-IE framework (see also Section~\ref{sec:npsem-dag}), the proof of the backdoor criterion \citep[][Theorem~3.2.2 and \S~11.3.3]{pearl2009causality} relies only on the g-formula~\eqref{eqn:truncated-distribution}, and therefore holds under the causal DAG model.
\end{remark}

The backdoor criterion can be operationalized using the following simple rules \citep[\S~11.3.1]{pearl2009causality}, such that after conditioning on \( L \),
\medskip 

\begin{enumerate}
    \item[{\bf Rule 1:}] All non-causal  paths between \( A \) and \( Y \) are blocked, and no new ones are introduced.\\[-8pt]
    \item[{\bf Rule 2:}] All causal  paths from \( A \) to \( Y \) remain open (unblocked).
\end{enumerate}
\medskip

\begin{figure}[ht]
\centering

% Subfigure (a)
\begin{subfigure}[t]{0.45\textwidth}
\centering
\begin{tikzpicture}[>=Stealth, node distance=1cm, every node/.style={draw, circle}]
  \node (X)                         {\(A\)};
  \node (S) [right=of X]            {\(S\)};
  \node (Y) [right=of S]            {\(Y\)};
  \node (Z) [below=of S]            {\(Z\)};
    % Dashed circle for epsilon_S
  \node[draw, circle] (eS) [above=of S] {\(\varepsilon_S\)};

  \draw[->] (X) -- (S);
  \draw[->] (S) -- (Y);
  \draw[->] (S) -- (Z);
% Dashed arrow from epsilon_S to S
  \draw[->] (eS) -- (S);
\end{tikzpicture}
\caption{A DAG with post-treatment variables \( S \),  \( Z \), and $Y$.}
\end{subfigure}
\hfill
% Subfigure (b)
\begin{subfigure}[t]{0.45\textwidth}
\centering
\begin{tikzpicture}[>=Stealth, node distance=1cm, every node/.style={draw, circle}]
  \node (X)                         {\(A\)};
  \node (Y) [right=2cm of X]        {\(Y\)};
  \node (Z) [below=1cm of X, xshift=1cm] {\(Z\)};
  \node[draw, circle] (eS) [above=1cm of X, xshift=1cm] {\(\varepsilon_S\)};
  \draw[->] (X) -- (Y);
  \draw[->] (X) -- (Z);
  \draw[->] (eS) -- (Y);
  \draw[->] (eS) -- (Z);
\end{tikzpicture}
\caption{The same DAG after marginalizing over \( S \).}
\end{subfigure}

\caption{Left: the original model including the mediator \( S \). Right: the simplified model where \( S \) is removed and its influence is captured via direct arrows.}
\label{fig:dag}
\end{figure}
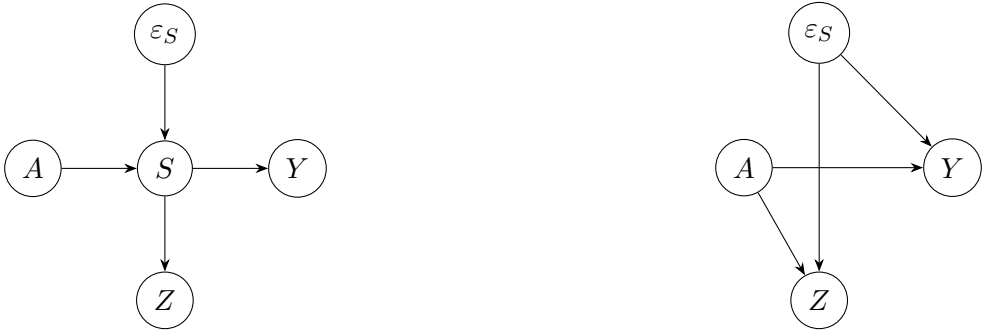

The first rule of the backdoor criterion suggests that adjusting for post-treatment variables may introduce bias. To illustrate, consider the causal DAG in Figure~\ref{fig:dag}(a), 
where the treatment \(A\) has two descendants, \(S\) and \(Z\), 
in addition to the primary variables $A$ and $Y$.  Adjusting for \( S \) clearly violates Rule 2 above. Now consider adjusting for \( Z \). On the surface, \( Z \) does not lie on any path from \( A \) to \( Y \), and \( S \) is not a collider on the path \( A \rightarrow S \rightarrow Y \). However, every observed variable has an associated error term, and in panel (a) we 
explicitly display \(\varepsilon_S\), the error for \(S\), since it plays a key role here. 
Conditioning on \(Z\) induces a dependence between \(A\) and \(\varepsilon_S\), 
which in turn opens the non-causal path 
\(A \leftrightarrow \varepsilon_S \to S \to Y\). 
Equivalently, this corresponds to the undirected path 
\(A - \varepsilon_S - S - Y\) in the induced moral graph. 
Another way to see this is by marginalizing out \( S \) and examining the resulting causal DAG in Figure~\ref{fig:dag}(b). In this graph, it becomes clear that \( Z \) is a collider on the non-causal path \( A \rightarrow Z \leftarrow \varepsilon_S \rightarrow S \), and hence adjusting for \( Z \) may introduce collider bias into the causal effect estimate \citep[\S~11.3.1]{pearl2009causality}.

% \tblue{Add results that it always exists.}

When the backdoor criterion fails to apply directly, identification may still be possible through Pearl's \emph{do-calculus} \citep{pearl1995causal}, a set of three transformation rules that enable the conversion of interventional expressions (e.g., \( P_{\text{do}(a)}(Y)  \)) into estimable quantities from the observational distribution. While the backdoor criterion suffices for simple DAGs, complete graphical criteria for covariate adjustment in ancestral and partially directed graphs have been developed \citep{Perkovic2018CompleteAdjustment}. In fact, any identifiable interventional distribution under a causal DAG model that may involve unmeasured nodes can be derived through recursive application of the do-calculus rules, as formalized by the ID algorithm \citep{tian2002identification,huang2006pearl,shpitser2006identification}. See \cite{tian2010identifying} for a review of these results.

% \tblue{Double check the ID algorithm assumes the causal DAG, not NPSEM-IE. Also the backdoor criterion. The Pearl 1995 paper assumes NPSEM-IE; as well as Pearl 2009 p.~85 Theorem 3.4.1. }

It is important to note that this result pertains to identifiability under the causal DAG model alone. Given additional assumptions, identification may still be possible even when the ID algorithm fails to apply. For instance, one can show that the ACE is not identifiable under the instrumental variable model alone \citep{balke1997nonparametric}. However, it can still be identified under additional effect homogeneity assumptions \citep[][Theorem 1]{wang2018bounded}; see also \citet[][Theorem 1]{hartwig2023average} and \citet[][Theorem 1]{dong2025marginal}.

\begin{remark}
\label{remark:id}
  There are more sophisticated algorithms that also identify counterfactual quantities, such as those discussed in Remark~\ref{remark:counterfactual} \citep{shpitser2007what,shpitser2008complete}. These algorithms, unlike the ID algorithm mentioned above, rely on the NPSEM-IE framework.
\end{remark}

% \subsection*{Conclusion}

% The DAG framework offers a robust language for causal reasoning and a toolkit—including d-separation, backdoor/frontdoor criteria, and do-calculus—for identifying causal effects. It allows analysts to represent assumptions visually, derive testable implications, and assess the identifiability of interventions. Crucially, DAG-based inference separates the causal logic from statistical estimation, ensuring that identification is grounded in substantive assumptions about the data-generating process.

% \bigskip
% \bigskip
% \bigskip

\section{Translating Between the Three Causal Frameworks}
\label{sec:translating}

The question of which causal framework is best suited for practical applications has been the subject of extensive debate in the causal inference literature \citep[e.g.][]{pearl1995causal, rubin2004direct, lauritzen2004discussion,robins2011alternative,richardson2023potential}; see also \citet[][p.~106]{pearl2009causality}. While each approach, including potential outcomes, structural causal models,  and graphical models, offers unique strengths, efforts have been made to reconcile and translate among them. Much like natural languages, certain causal concepts are more naturally expressed in one framework than another, making such translations both valuable when possible but, at times, inherently limited. In this section, we review key attempts to bridge these perspectives. Detailed comparisons of these frameworks are deferred to Section \ref{sec:comparison}.

\subsection{The Potential Outcomes and NPSEM-IE Frameworks}
  
% There exists a deep and formal connection between the potential outcomes and the NPSEM-IE frameworks. Both frameworks provide semantics for reasoning about counterfactuals, but they approach the problem from different perspectives. The potential outcomes framework defines counterfactuals as primitives indexed by hypothetical interventions, while NPSEM-IEs adopt a structural, functional view of how variables are generated. 

\paragraph{From NPSEM-IE to Potential Outcomes.} 
Under an NPSEM-IE, each endogenous variable \( V_j \) is determined by a deterministic function of its parents and an exogenous noise variable \( \varepsilon_j \), with the collection \( \{ \varepsilon_j \} \) assumed to be jointly independent. This structural system implicitly defines a rich set of potential outcomes.  For example, consider the NPSEM in \eqref{eq:npsem}. In this case, the potential outcome \( Y(l,a) \) is obtained by intervening to set \( L = l \), \( A = a \),  and evaluating the outcome function at the realized value of the noise term \( \varepsilon_Y \) \citep[][\S3.6.3, \S 7.1]{pearl2009causality}:
\begin{equation}
    \label{eqn:npsem-po1}
    Y(l,a) = f_Y(l, a, \varepsilon_Y);
\end{equation}
see also \cite{strotz1960recursive}.
In fact, researchers have proposed expressing NPSEMs using potential outcomes notation to clearly distinguish them from standard statistical regression models \citep[e.g.][]{imbens:2014, richardson:robins:2014, richardson2023potential}.
% \tblue{todo: check this reference}
 For example, the NPSEM in \eqref{eq:npsem} can be represented as
\begin{align}
L &= f_L(\varepsilon_L), \notag \\
A(l) &= f_A(l, \varepsilon_A), \label{eq:rewrite}  \\
Y(l,a) &= f_Y(l, a, \varepsilon_Y), \notag
\end{align}
where the potential outcomes notation explicitly emphasizes the causal interpretation of the structural equations.  

More compactly, we can also define the marginal potential outcome by absorbing \(L\) into its realized value under the intervention:
\begin{equation}
    \label{eqn:npsem-po2}
    Y(a) \;=\; Y(L,a) \;=\; f_Y(L, a, \varepsilon_Y).
\end{equation}

One can also recursively define nested potential outcomes from NPSEMs. For example, consider a simple model with a treatment \( A \), a mediator \( M \), and an outcome \( Y \), governed by:
\begin{align*}
    A &= f_A(\varepsilon_A), \\
    M &= f_M(A, \varepsilon_M), \\
    Y &= f_Y(A, M, \varepsilon_Y).
\end{align*}
 In this setting, the nested counterfactual \( Y(a, M(a')) \) represents the value of \( Y \) if treatment were set to \( a \), and the mediator were assigned the value it would attain under treatment \( a' \). This can be expressed as:
\begin{equation}
    \label{eqn:npsem-po3}
Y(a, M(a')) = f_Y\big(a, f_M(a', \varepsilon_M), \varepsilon_Y\big).
\end{equation}

\paragraph{Implications of Missing Input Variables in NPSEM-IE}  In the NPSEM defined in \eqref{eqn:npsem-ie}, the set of input variables \( U_j \) for a given node \( V_j \) often does not include all variables that precede \( V_j \) in the causal ordering. The omission of these variables reflects an implicit \emph{exclusion restriction} on the corresponding potential outcomes \citep[][p.~101, eqn.~(3.55)]{pearl2009causality}:
\[
V_j(u_j) = V_j(u_j, s),
\]
where \( S \subset V \) is disjoint from \( U_j \). This equality states that setting variables in \( S \) does not affect the value of \( V_j \) once the variables in \( U_j \) are held fixed, thereby encoding a form of causal irrelevance.

As a concrete example, consider a simple instrumental variable (IV) model with an unmeasured confounder \( U \) that affects both the treatment \( A \) and the outcome \( Y \). Let \( Z \) be an observed instrument that influences \( A \) but has no direct effect on \( Y \). The structural equations may be specified as:
\[
\begin{aligned}
Z &= \varepsilon_Z, \\
U &= \varepsilon_U, \\
A &= f_A(Z, U, \varepsilon_A), \\
Y &= f_Y(A, U, \varepsilon_Y).
\end{aligned}
\]
In this setup, \( U \) introduces unmeasured confounding between \( A \) and \( Y \). The exclusion restriction required for IV identification is encoded in the structural equation for \( Y \): the variable \( Z \) is omitted from the input set in the SEM for $Y$. This implies that, for any \( a \) and \( z \),
\[
Y(z, a, U) = Y(a, U),
\]
or, equivalently,
\begin{equation}
    \label{eqn:exclusion}
    Y(z,a) = Y(a).
\end{equation}
Eqn. \eqref{eqn:exclusion} is often known as the exclusion restriction in the IV literature \citep{angrist1996identification}.

\paragraph{From Potential Outcomes to NPSEM-IE via one-step-ahead counterfactuals} 
\citet{robins2011alternative, richardson2013single} formalize the converse direction by showing how one can construct an NPSEM-IE from a given system of potential outcomes using what they term \emph{one-step-ahead counterfactuals}, i.e., the potential outcome of a variable under interventions on all of its parents \citep{robins1986new}.  This ``counterfactual interpretation of the error terms in structural equation models'' was also noted in \citet[][\S 7.4.5]{pearl2009causality} and referred to as the ``canonical representation of structural causal models'' by \citet[][\S 3.4]{peters2017elements}.

To illustrate, consider a binary covariate \( L \) and a binary treatment variable \( A \). Suppose the available potential outcomes include \( A(l) \) and \( Y( l,a) \) for \( l \in \{0,1\} \) and \( a \in \{0,1\} \). One can reconstruct the NPSEM as follows: 
\begin{alignat}{2}
\varepsilon_L &= L,             &\qquad f_L(\varepsilon_L) &= \varepsilon_L, \notag \\
\varepsilon_A &= (A(1), A(0)),  &\qquad f_A(l, \varepsilon_A) &= 
\begin{cases}
    A(1), & \text{if } l = 1, \\
    A(0), & \text{if } l = 0,
\end{cases} \label{eqn:npsem-canonical} \\
\varepsilon_Y &= \big(Y(1,1), Y(1,0), Y(0,1), Y(0,0)\big), 
&\qquad 
f_Y(a,l,\varepsilon_Y) &= 
\begin{cases}
    Y(1,1), & \text{if } l=1,\, a=1, \\
    Y(1,0), & \text{if } l=1,\, a=0, \\
    Y(0,1), & \text{if } l=0,\, a=1, \\
    Y(0,0), & \text{if } l=0,\, a=0.
\end{cases} \notag
\end{alignat}
In this construction, each exogenous variable encapsulates the relevant potential outcomes for the corresponding endogenous variable, while the structural functions act as selectors based on the realized values of upstream variables.

% \paragraph{Comparing the NPSEM-IE with the Potential Outcomes.}   

% \paragraph{Summary.} 
% These two frameworks are therefore mathematically equivalent under appropriate assumptions, with NPSEM-IEs providing a functional representation that generates potential outcomes, and potential outcomes defining mappings from which an NPSEM-IE can be reconstructed. This equivalence forms the foundation for many developments in modern causal inference, including identification theory and counterfactual reasoning.

% \begin{figure}[ht]
% \centering
% \begin{tikzpicture}[->, node distance=2cm, thick]
%   \node[draw, circle] (A) {A};
%   \node[draw, circle, right of=A] (M) {M};
%   \node[draw, circle, right of=M] (Y) {Y};

%   \node[below of=A, yshift=0.7cm] (FA) {\small $A = f_A(\varepsilon_A)$};
%   \node[below of=M, yshift=0.7cm] (FM) {\small $M = f_M(A, \varepsilon_M)$};
%   \node[below of=Y, yshift=0.7cm] (FY) {\small $Y = f_Y(A, M, \varepsilon_Y)$};

%   \draw[->] (A) -- (M);
%   \draw[->] (M) -- (Y);
%   \draw[->] (A) to[bend right=15] (Y);
% \end{tikzpicture}
% \caption{Counterfactuals in NPSEM-IE are constructed by recursively substituting values (e.g., $A = a$) into structural equations using shared exogenous errors.}
% \label{fig:npsem_counterfactuals}
% \end{figure}

\subsection{NPSEM-IE  and  Causal DAG}
\label{sec:npsem-dag}

\subsubsection*{From an NPSEM-IE to a Causal DAG}
\label{sec:npsemtodag}

Consider the NPSEM \eqref{eqn:npsem-ie}.
To construct a graph from this NPSEM, one proceeds as follows:  
(1) Draw a node for each endogenous variable \( V_j \);  
(2) For each variable \( V_j \), draw directed edges from every variable that appears as an argument in \( f_j \) to \( V_j \). These variables are the parents of \( V_j \) in the graph.  
If the resulting graph is acyclic, i.e., it forms a DAG, then the corresponding NPSEM is referred to as an \emph{acyclic NPSEM}. In other words, any acyclic NPSEM induces a DAG in a purely graphical sense. The following proposition establishes the connection between the independent error and autonomy assumptions in an NPSEM-IE and the causal DAG model.

\begin{proposition}
\label{prop:npsem-dag}
Any NPSEM-IE implies a unique causal DAG model. In particular, it induces  
(1) an observational distribution that factorizes according to the DAG induced by the NPSEM;  
(2) interventional distributions that satisfy the truncated factorization formula \eqref{eqn:truncated-distribution}.%
\end{proposition}

The first claim in Proposition~\ref{prop:npsem-dag} follows from \citet[][Theorem 1.4.1]{pearl2009causality}; see also \citet[][Proposition 6.31]{peters2017elements}. The second claim holds because, under the autonomy assumption, the structural equation for the intervened variable is replaced while all others remain unchanged; see also \citet[][\S2]{strotz1960recursive}, \citet[][p.~72 eqn. (3.10)]{pearl2009causality} and \citet[][p.~109, eqn. (6.7)]{peters2017elements}.

\begin{remark}
    In fact, a stronger claim holds: 
the Single-World Intervention Graph (SWIG) model (see \S\ref{section:swig} below) with interventions on all variables also implies a causal DAG model \citep[][Theorem 11]{richardson2023potential}. This is stronger because 
the NPSEM-IE is a strict submodel of the SWIG model associated with the same DAG. The converse is not necessarily true. For example, the SWIG in Figure \ref{fig:swig} implies the so-called extended g-formula \citep{robins2004effects} that
\begin{equation}
    \label{eqn:extended-g-formula}
    P(A=a^*, L=l, Y(a) = y) = P(L=l) P(A=a^*\mid L=l) P(Y=y\mid A=a, L=l) 
\end{equation}
for any values of $a^*, a, l,$ and $y$. 
Although this implies the truncated factorization formula \eqref{eqn:truncated}, the converse does not hold.
\end{remark}

\begin{remark}
     The extended g-formula \eqref{eqn:extended-g-formula} can be used to identify counterfactual quantities such as the  ETT. As noted later in Section \ref{sec:single-vs-cross-world}, the ETT is not  expressible in the causal DAG framework.
\end{remark}

% \tblue{Another way to see this is through SWIG: (1) Thomas' slide 57 says the Independences arising from a SWIG imply all of the
% identification results that hold in the do-calculus of Pearl (1995);
% see also Spirtes et al. (1993), which seems to suggest that SWIG implies causal DAG; (2) I also remember that the SWIG paper says that NPSEM-IE implies the SWIG.}

\subsubsection*{From Causal DAGs to NPSEM-IE}

Let $\mathcal{G}$ be a causal DAG over variables $V = \{V_1, \dots, V_p\}$, and let $P(V_1, \dots, V_p)$ be a joint distribution that satisfies the Markov factorization with respect to $\mathcal{G}$:
\[
P(V_1, \dots, V_p) = \prod_{j=1}^p P(V_j \mid \text{Pa}_{\mathcal{G}}(V_j)).
\]
The Functional Representation Lemma \citep[][p.~626, Appendix B]{el2011network} ensures that any random variable with a known conditional distribution can be expressed as a measurable function of its parents and an independent noise variable. A concrete construction is provided in \citet[][Appendix C.9]{peters2017elements}.

% \begin{lemma}[Functional Representation Lemma]
% Let $(X, Y)$ be random variables on a probability space, where $X \in \mathcal{X}$ and $Y \in \mathcal{Y}$, and suppose the conditional distribution $P_{Y \mid X}$ is known. Then there exists:
% \begin{enumerate}
%     \item A measurable function $f: \mathcal{X} \times [0,1] \to \mathcal{Y}$,
%     \item A random variable $\varepsilon \sim \text{Uniform}[0,1]$, independent of $X$,
% \end{enumerate}
% such that $Y = f(X, \varepsilon)$.
% \end{lemma}

Applying the functional representation lemma    recursively to each conditional distribution $P(V_j \mid \text{Pa}_{\mathcal{G}}(V_j))$, we obtain an NPSEM:
\begin{equation}
    \label{eqn:npsem-ie-2}
    V_j = f_j(\text{Pa}_{\mathcal{G}}(V_j), \varepsilon_j), \quad j = 1, \dots, p,
\end{equation}
where each \( f_j \) is a measurable function, each \( \varepsilon_j \sim \text{Uniform}[0,1] \), and the noise variables \( \varepsilon_1, \dots, \varepsilon_p \) are jointly independent. The NPSEM \eqref{eqn:npsem-ie-2} with independent errors is compatible with both the causal DAG \( \mathcal{G} \) and the joint distribution \( P \). A similar construction for the discrete case can be found in \citet{druzdzel1993causality}.

\begin{proposition}
    Consider a causal DAG model over random variables $V = \{V_1, \ldots, V_p\}$ as in Definition \ref{def:causal-dag}. Then there exists an NPSEM-IE that implies the causal DAG model. 
\end{proposition}

\begin{remark}
    Although this construction is always possible under the assumptions of the functional representation lemma, the resulting NPSEM-IE is generally not unique, as different $(f_j, \varepsilon_j)$ pairs may induce identical observational distributions and factorizations.
\end{remark}

% \subsection*{Comparing the NPSEM-IE with Causal DAGs}

\subsection{Potential Outcomes and DAGs: Single World Intervention Graphs (SWIGs)}
\label{section:swig}

The potential outcomes and DAG frameworks can be unified through \emph{single world intervention graphs} (SWIGs) \citep{richardson2013single}. SWIGs explicitly encode potential outcomes on the graph by splitting nodes into fixed (intervened) and random components. This representation enables two key features: (1) conditional independence statements involving both observed and potential outcomes can be read using a modified d-separation rule, and (2) identification results equivalent to do-calculus can be derived from these conditional independence conditions via a simple reformulation known as the potential outcome (po) calculus \citep{malinsky-19b,shpitser2022multivariate}. This reformulation shows that from the perspective of potential outcomes, the po-calculus (and hence the do-calculus)  
is reducible to conditional independence, recursive substitution and causes only influence descendants in the SWIG. For completeness, we note that SWIGs can also be defined via a local Markov property 
with respect to a total ordering of the variables, without explicit reference 
to a graphical representation \citep[see][]{richardson2023potential}.
The same potential outcome model was originally termed the
Finest Fully Randomized Causally Interpretable Structured Tree Graph (FFRCISTG); see \citet{robins1986new} and \citet[][Appendix C]{richardson2013single}; note that the ``tree graphs'' referred to here are ``event trees,'' not DAGs.

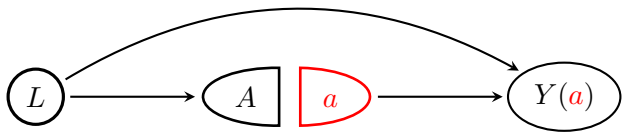
\begin{figure}[!htbp]
\centering
	\begin{tikzpicture}[>=stealth, line width = 1pt,   node distance=1.5cm,
	% \tikzstyle{format} = [draw, ultra thick, minimum size=6mm
	pre/.style={->,>=stealth,black, inner sep=8pt}]
\begin{scope}
\tikzset{swig vsplit={gap=8pt,
line color right=red}}
\node[name=A,pre, shape=swig vsplit]{
\nodepart{left}{$A$}
\nodepart{right}{\textcolor{red}{$a$} }}; 
%	\begin{scope}
%	\node[name=A,shape=ellipse splitb, ellipse splitb/colorleft=black, ellipse splitb/colorright=red, 
%	ellipse splitb/innerlinewidthright = 0pt,  %Remove the inner 'line'
%	/tikz/ellipse splitb/linewidthright = 1pt,   %Make Right line same width as left
%	ellipse splitb/gap=3pt, style={draw},rotate=90] {
%		{\rotatebox{-90}{$A$\;}}
%		\nodepart{lower}
%		{\rotatebox{-90}{\;\textcolor{red}{$a$}}}
%	};
	% \node[name=D,shape=ellipse splitb, below right=0.65cm and 1.8cm of Z,  ellipse splitb/colorleft=black, ellipse splitb/colorright=red, 
	% ellipse splitb/innerlinewidthright = 0pt,  %Remove the inner 'line'
	% /tikz/ellipse splitb/linewidthright = 1pt,   %Make Right line same width as left
	% ellipse splitb/gap=3pt, style={draw},rotate=90] {
		% {\rotatebox{-90}{$D(\tred{z})$\;}}
		% \nodepart{lower}
		% {\rotatebox{-90}{\;\textcolor{red}{$d$}}}
	% };
	 \node[thick, name=Y,shape=ellipse,style={draw}, right= 1.8cm of A, outer sep=0pt, text width = 8mm]
	 {$Y({\color{red}{a}})$
	 };
	\draw[pil,->] (A) to (Y);
	\node[est, left= 1.8cm of A] (L) {$L$};
        \draw[pil,->] (L) to (A);
        \draw[pil,->] (L) to [bend left] (Y);
	% \draw[pil, ->] (2.7,2)  to[bend right]  (-0.5,0.4);
	% \draw[pil, ->] (3.1,1.7)  --  (3.1,0.4);
	% \draw[pil, <->] (3.1,-0.5)  to[bend left]  (-0.5,-0.4);
	% \path[pil] (V) edgenode {} (Y);
	% \node[shade, below = of D] (U) {$U$};
	% \path[pil] (U) edgenode {} (3.3,-0.5);
	% \path[pil] (U) edgenode {} (Y);
	% %\node[below right = 8mm of anew2]{(a)};
	\end{scope}
	\end{tikzpicture}
    \caption{The single world intervention graph (template) corresponding to intervening on $A$ in the causal DAG in Figure \ref{fig:dag_complete}.}
\label{fig:swig}
\end{figure}

For example, consider intervening on the variable $A$ in the DAG presented in Figure~\ref{fig:dag_complete}. As illustrated in Figure~\ref{fig:swig}, a SWIG can be constructed in two steps. In the first step, one splits each intervention node, in this case, $A$, into two nodes: a random node $A$, which inherits all incoming arrows from the original DAG, and a fixed node $a$, which inherits all outgoing arrows from the original DAG. In the second step, all descendants of the fixed node $a$ are relabeled as potential outcomes, expressed as functions of the fixed node's value. Formally, the graph in Figure~\ref{fig:swig} is a \emph{template} (or a graph-valued function), since the fixed node $a$ can take multiple values, such as $a = 0$ or $a = 1$.

One can read off conditional independence relationships from a SWIG using the d-separation rule, with the extension that a d-connecting path must \emph{not} contain fixed (red) nodes. For example, from the SWIG in Figure \ref{fig:swig}, one can directly read off the weak ignorability condition \eqref{eqn:weak-ignorability} that $A\ind Y(a)\mid L, \text{for } a=0,1$ using the rule of d-separation.

Graphs like the one in Figure~\ref{fig:swig} are called \emph{single-world} because potential outcomes corresponding to \emph{different} hypothetical interventions on the same variable never appear in the same SWIG. Consequently, a SWIG does not encode assumptions about relationships between potential outcomes across different interventions on the same variable. For example, the SWIG in Figure~\ref{fig:swig} does \emph{not} imply the strong ignorability condition $A \ind (Y(1), Y(0)) \mid L$. We discuss the distinction between single-world and cross-world assumptions in more detail in Section~\ref{sec:single-vs-cross-world}.

\begin{figure}[htbp]
    \centering
    \begin{subfigure}[b]{0.42\textwidth}
        \centering
        \includegraphics[width=.7\linewidth]{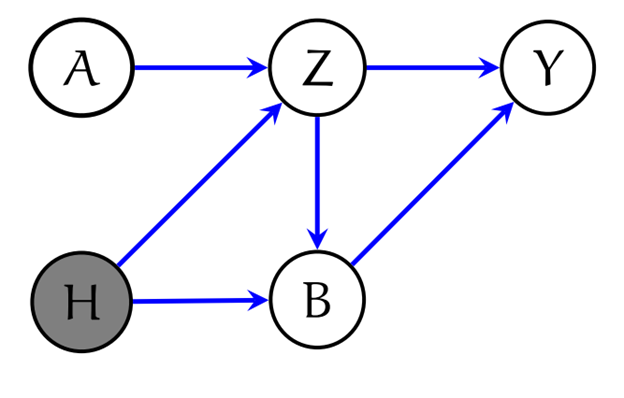}
        \caption{A complex causal DAG model.}
        \label{fig:sub1}
    \end{subfigure}
    \hfill
    \begin{subfigure}[b]{0.54\textwidth}
        \centering
        \includegraphics[width=.7\linewidth]{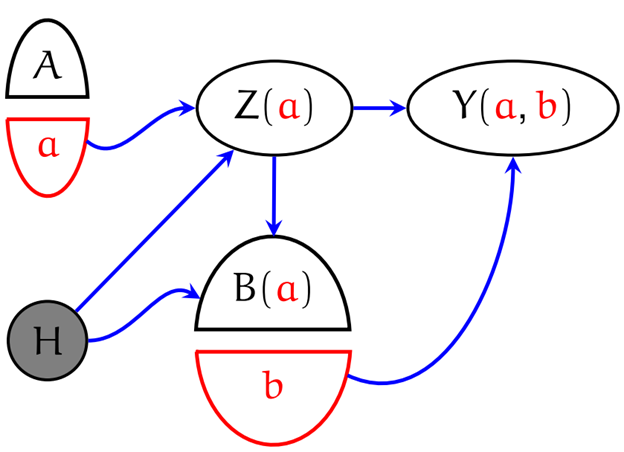}
        \caption{The corresponding SWIG model with interventions on $A$ and $B$.}
        \label{fig:sub2}
    \end{subfigure}
    \caption{Does the causal DAG model imply $Y(a,b) \ind B \mid Z, A=a$?}
    \label{fig:side_by_side}
\end{figure}

SWIGs are useful not only for validating well-known conditional independence assumptions that hold under the causal DAG model, but also as graphical tools for assessing conditional independence statements that are not obvious from a causal DAG model. For example, consider the conditional independence statement
\begin{equation}
    \label{eqn:cond-ind}
    Y(a,b) \ind B \mid Z, A=a.
\end{equation}
\citet[][Ex. 11.3.3]{pearl2009causality} claimed that the causal DAG model in Figure~\ref{fig:sub1} does not imply \eqref{eqn:cond-ind}. However, as can be seen from the SWIG in Figure~\ref{fig:sub2}, the relationship $Y(a,b) \ind B(a) \mid Z(a), A=a$ holds, which implies \eqref{eqn:cond-ind} through the consistency assumption.

% todo: \linbo{Are SWIGs the same as the causal DAG model? What is their relationship? Ask for tikz code from Thomas.}

% Add a remark on single versus cross world assumptions. 

% Pearl's ‘do' notation is not as rich as the
% language of potential outcomes as there is no way to directly express $P(Y (a)\mid A =a^\prime)$ in
% terms of his ‘do' notation \citep{shpitser2022multivariate}. Double check this with Thomas: Do-notation is single-world. 

\section{Comparison Between the Three Frameworks}
\label{sec:comparison}

\subsection{Minimalist vs. Mechanistic Perspectives on Causality}

Although, as we have shown, these frameworks often lead to overlapping mathematical results, they arise from different philosophical traditions and emphasize different aspects of the causal enterprise \citep[\S~1.4]{pearl2009causality}.  
The potential outcomes framework adopts a \emph{minimalist} or design-based view of causality, focusing on comparisons between hypothetical outcomes under different treatment conditions, such as \( Y(1) - Y(0) \). This approach does not require specifying how all outcomes are generated from other variables in the system. It is particularly attractive in settings where causal identification relies on features of the study design, such as randomization or the use of instrumental variables, rather than on a full model of the data-generating process. Because of its minimal assumptions and broad applicability, this framework has been described by some authors as the most general approach to causal inference \citep{imbens2015causal,rosenbaum2010design}.

In contrast, NPSEM-IEs adopt a more \emph{mechanistic} perspective, rooted in classical physics and engineering, where causality is represented through deterministic functional relationships between variables and exogenous noise terms. In this tradition, causal systems are fully specified via structural equations, and randomness reflects ignorance about latent factors rather than intrinsic stochasticity. This view aligns with Laplace's deterministic vision of natural laws \citep{laplace1825essai} and provides a rich language for reasoning about complex causal systems.

Causal DAGs provide a bridge between this minimalist perspective and mechanistic modeling. By encoding assumptions about independencies, DAGs help justify potential-outcomes-based identifiability conditions (e.g., ignorability) while avoiding over-specification of functional forms. For instance, a DAG makes it transparent why randomization implies that no arrows point into $A$. This, in turn, ensures the independence $A \ind Y(a), \, a \in \{0,1\}$, which can be read off from a SWIG intervening on $A$. Moreover,
DAGs clarify backdoor adjustment sets for observational studies \citep{pearl2009causality}.  In contrast to NPSEMs, all relationships in causal DAG models are formulated as inherently stochastic. As \citet[p. 26]{pearl2009causality} observed, this reflects the modern conception of physics, sometimes described as quantum mechanical, in which probabilistic laws replace strict determinism as a more realistic description of physical processes. At the same time, DAGs retain the ability to articulate complex structural relationships, making them a powerful tool for both identification and hypothesis generation.

As a caveat, we note that although Pearl’s analogy draws on the idea that modern physics views natural laws as fundamentally probabilistic rather than strictly deterministic, it should be interpreted only metaphorically: causal DAGs are classical probabilistic models that admit hidden variable representations, whereas quantum systems in general do not. Quantum phenomena require richer nested probabilistic frameworks, sometimes viewed as supermodels of the classical theory, that extend beyond the scope of standard DAGs \citep{richardson2023nested}.

% causal DAG models are formulated directly in probabilistic terms and therefore align with the modern, quantum-mechanical perspective, according to which fundamental physical laws are inherently stochastic. 

% \tblue{Talk about the causal inference roadmap here.}

% \tblue{Read Peters section 6.9.2.}

In what follows, we compare these three frameworks in greater depth, highlighting their respective assumptions, strengths, and areas of applicability.

\subsection{Single-World vs. Cross-World Causal Inference}
\label{sec:single-vs-cross-world}

\paragraph{Difference in Expressive Power for Counterfactual Quantities.} The key difference between the causal DAG framework with its associated 
``do'' notation and the other two frameworks (NPSEM-IE and potential outcomes) 
is that the former does not provide a way to directly express counterfactual 
quantities that condition on observed variables rather than interventions, 
such as $P(Y(a) \mid A = a')$ \citep[cf.][§1.4.4]{pearl2009causality}. 
An important example is the ETT,
$\mathbb{E}[Y(1)-Y(0)\mid A=1]$, which is not expressible in the ``do'' notation.
Compared with the causal DAG model, the NPSEM-IE model also imposes many more assumptions due to the additional cross-world assumptions.  We now elaborate this point.

\paragraph{Implications of the Independent Error Assumptions in NPSEM-IE}

Because the independent error assumptions in NPSEM, when taken at face value, resemble the independent error assumptions commonly made in statistical modeling, they are often perceived as mild by statisticians who are not familiar with the causal modeling framework. Interestingly, many of the same statisticians would view the assumption of no unmeasured confounding as an extremely strong assumption. Before the recent evolution in causal inference, it was standard practice in statistics to avoid causal claims altogether, precisely because of the concern that confounding renders correlation insufficient for establishing causation \citep[][Preface]{pearl2016causal}.

However, as we demonstrate below, the independent error assumptions in an NPSEM-IE are, in fact, much stronger than the assumption of no unmeasured confounding as in eqn. \eqref{eqn:weak-ignorability}.

Consider the NPSEM-IE in eqn. \eqref{eq:npsem}. The independent error assumption states that \( \varepsilon_L \ind \varepsilon_A \ind \varepsilon_Y \). Now consider the canonical representation in eqn. \eqref{eqn:npsem-canonical}. Under this construction, the  joint independence of the exogenous variables implies
\begin{equation}
    \label{eqn:ie}
    L \ind \left(A(l) : l \in \mathcal{L}\right) \ind \left(Y(l,a) : a \in \{0,1\}, l \in \mathcal{L}\right),
\end{equation}
where \( \mathcal{L} \) denotes the support of \( L \) \citep[][p.~101 footnote 14]{pearl2009causality}. Basic algebra then yields the strong ignorability condition,
\begin{equation}
    \label{eqn:strong-ignorability}
    A \ind (Y(1), Y(0)) \mid L,
\end{equation}
which is strictly stronger than the weaker no-unmeasured-confounding assumption in eqn. \eqref{eqn:weak-ignorability}.

In fact, we have the following general result due to \citet[][p.~101, eqn. (3.56) \& footnote 14]{pearl2009causality}.
\begin{proposition}
\label{prop:npsem-ie}
  In NPSEM \eqref{eqn:npsem-ie}, joint independence of the $p$ exogenous errors 
$\varepsilon_j, j=1,\ldots,p$ is equivalent to joint independence of the following $p$ collections of potential outcomes 
$\{V_j(\text{pa}(V_j)=a) : a \in \mathcal{A}_{\text{pa}(V_j)}\}, \ j=1,\ldots,p$, where $\mathcal{A}_{\text{pa}(V_j)}$ denotes the 
set of all possible assignments to the parents of $V_j$.
\end{proposition}

Note that joint independence of the exogenous errors does 
{\em not} imply that the variables $V_j(\text{pa}(V_j)=a)$,
$a \in \mathcal{A}_{\text{pa}(V_j)}$ are independent of one another; thus for example this does not imply $Y(a=0) \ind 
Y(a=1)$.

To understand the number of extra assumptions imposed by the NPSEM-IE, consider again the NPSEM-IE in eqn.~\eqref{eq:npsem} and a simple case where \( L \), \( A \), and \( Y \) are all binary. A full model for the joint distribution of \( (L, A(l), Y(l,a)) \), for \( a \in \{0,1\} \) and \( l \in \{0,1\} \), has dimension \( 2^{1 + 2 + 4} - 1 = 127 \). The no-unmeasured-confounding assumption in eqn. \eqref{eqn:weak-ignorability} reduces this to a 123-dimensional model. In contrast, the strong ignorability assumption in eqn. \eqref{eqn:strong-ignorability} further reduces the model to 121 dimensions, and the independent error assumption of the NPSEM-IE implies a model of only 19 dimensions!
In other words, NPSEM-IE imposes an additional \( 123 - 19 = 104 \) constraints, all of which are not necessary for the identification of the ACE \citep{robins2011alternative}. Most of these constraints, such as \( A(l=1) \ind Y(l=0,a) \), are examples of \emph{cross-world assumptions}, which remain untestable even under randomized treatment assignment.  These assumptions violate the principle of ``no causation without manipulation'' \citep{holland1986statistics} and undermine the connection between experimentation and causal inference \citep{robins2003semantics,robins2011alternative}. 
% See Section~\ref{section:swig} for further discussion of cross-world assumptions.

\paragraph{When the NPSEM-IE Framework Leads to Stronger Identification Results} 

In some cases, however, the formulation and additional assumptions in NPSEM-IE indeed lead to stronger identification results. Here, we provide several examples, each representing a different scenario in which these stronger results hold under the NPSEM-IE but not under alternatives such as the causal DAG model.

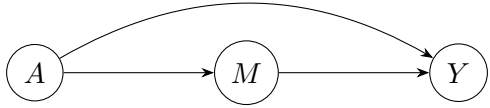
\begin{figure}[ht]
\centering

% Subfigure (a): Mediation
% \begin{subfigure}[t]{0.45\textwidth}
\centering
\begin{tikzpicture}[>=Stealth, node distance=2cm, every node/.style={draw, circle}]
  \node (A)                    {\(A\)};
  \node (M) [right=of A]       {\(M\)};
  \node (Y) [right=of M]       {\(Y\)};

  \draw[->] (A) -- (M);
  \draw[->] (M) -- (Y);
  \draw[->] (A) to[bend left=30] (Y);
\end{tikzpicture}
% \caption{Mediation model: \( A \to M \to Y \)}
% \end{subfigure}
% \hfill
% % Subfigure (b): IV model
% \begin{subfigure}[t]{0.45\textwidth}
% \centering
% \begin{tikzpicture}[>=Stealth, every node/.style={draw, circle}]
%   \node (Z)                    at (0, 0) {\(Z\)};
%   \node (A)                    at (2.5, 0) {\(A\)};
%   \node (Y)                    at (5, 0) {\(Y\)};
%   \node (U)                    at (3.75, 1.5) {\(U\)};

%   \draw[->] (Z) -- (A);
%   \draw[->] (A) -- (Y);
%   \draw[->] (U) -- (A);
%   \draw[->] (U) -- (Y);
% \end{tikzpicture}
% \caption{IV model with unmeasured confounder \( U \) between \( A \) and \( Y \)}
% \end{subfigure}

% \vspace{1em}

% \caption{Left: Mediation model used to derive the mediation formula. Right: IV model illustrating confounding and instrument exclusion.}
\caption{Mediation model used to derive the mediation formula. }
\label{fig:mediation-iv}
\end{figure}

\paragraph{Example 1: Mediation Formula.} In some cases, the additional independence assumptions in NPSEM-IE lead to new identification results. For example, 
Figure \ref{fig:mediation-iv} shows the mediation setting with $A$ the treatment, $M$ the mediator, and $Y$ the outcome of interest. An estimand of interest in this setting is known as the natural (or pure) direct effect \citep{robins:greenland:1992}, defined as $NDE = E[Y(1,M(0)) - Y(0)]$. The NDE can be identified under the NPSEM-IE using the mediation formula \citep{pearl2001direct}  
\[
\mathbb{E}[Y(1, M(0))] = \sum_m \mathbb{E}[Y \mid A = 1, M = m] \cdot P(M = m \mid A = 0).
\]
Specifically, it relies on the cross-world assumption that $M(0) \ind Y(1,m)$, which holds under the NPSEM-IE but not the causal DAG or SWIG assumptions. 

\paragraph{Example 2: Monotonicity Assumption.}  Some assumptions, such as the monotonicity assumption that \( Y(a = 0) \leq Y(a = 1) \) almost surely, cannot be expressed directly using a stochastic model like a causal DAG. Assumptions of this kind often lead to stronger identification results. For example, under the no unmeasured confounding assumption~\eqref{eqn:weak-ignorability} and the monotonicity assumption, the probability of necessity and sufficiency is identified as
\[
    PNS = P(Y = 1 \mid A = 1) - P(Y = 1 \mid A = 0).
\]
Other examples include the identification of the local average treatment effect in instrumental variable models~\citep{angrist1996identification}, tighter bounds on natural indirect effects in mediation analysis without assuming cross-world independencies~\citep{imai2010identification}, and the identification of, or tighter bounds on, the survivor average treatment effect~\citep{ding2011identifiability,  wang2017causal, wang2017identification}.

\paragraph{Example 3: Additive Noise Models.} 
\label{example:3}
In addition to strengthening identification of causal effects, NPSEM-IE assumptions can also aid in identifying the causal graph structure itself. This connects to our earlier discussion of Markov equivalence in Remark \ref{remark:markov}, where multiple DAGs may encode the same set of conditional independencies. In many cases, however, restricting the function class in NPSEM-IE resolves this ambiguity and leads to identifiability of the underlying causal structure. These assumptions are more naturally expressed in the language of NPSEMs. One such assumption is that the error terms are additive. Under this condition, many results have been established that lead to identification of the causal system in different scenarios. Below, we present several examples assuming an additive model for~\eqref{eqn:npsem-ie}:
\begin{equation}
    \label{eqn:npsem-additive}
    V_j = f_j(U_j) + \varepsilon_j,
\end{equation}
where \( U_j \) are the parents of \( V_j \), and \( \varepsilon_j \) are jointly independent.

\begin{enumerate}
    \item \textbf{(Linear, non-Gaussian)} \cite{shimizu2006linear} showed that if the functions \( f_j \) in the NPSEM-IE~\eqref{eqn:npsem-additive} are linear in \( U_j \) and \( \varepsilon_j \) with nonzero coefficients, and if the error terms \( \varepsilon_j \) are non-Gaussian with strictly positive density, then the underlying graph structure is identifiable from the observed data distribution. This model is known as the \emph{linear non-Gaussian acyclic model (LiNGAM)}.

    \item \textbf{(Nonlinear)} Consider the bivariate causal model
    \begin{equation}
        \label{eqn:bivariate}
        Y = f_Y(X) + \varepsilon_Y, \quad X \ind \varepsilon_Y.
    \end{equation}
   \cite{hoyer2008nonlinear} showed that, except in ``rare'' cases, a distribution does not admit an additive noise model structure in both directions simultaneously; see also \citet[][Theorem 4.5]{peters2017elements}.  In fact, if the noise \( \varepsilon_Y \) is Gaussian, then only \emph{linear} functions allow for additive noise models in both directions~\citep[][Corollary 1]{hoyer2008nonlinear}. More broadly, \cite{peters2014causal} showed that if the functions \( f_j \) in the NPSEM-IE~\eqref{eqn:npsem-ie} are three times differentiable and nonlinear in any component of \( U_j \)~\citep[see][Theorem 7.7 for a precise definition]{peters2014causal}, then the corresponding graph structure is identifiable from the observed data distribution. These results extend to the so-called \emph{post-nonlinear models} of the form \( V_j = g_j(f_j(U_j) + \varepsilon_j) \). These models also generally exist in at most one direction, except in rare cases~\citep{zhang2009postnonlinear}.
\end{enumerate}

Identification of the causal structure in linear models with Gaussian additive noise is more challenging. For example, consider again the bivariate causal model in~\eqref{eqn:bivariate}. If \( f_Y(X) \) is linear in \( X \), then there exist \( \beta \) and \( \varepsilon_X \) such that \( X = \beta Y + \varepsilon_X \) and \( Y \ind \varepsilon_X \) if and only if both \( \varepsilon_Y \) and \( X \) are Gaussian~\citep[][Theorem 4.2]{peters2017elements}. More generally, identification of the causal graph under linear models with Gaussian additive noise is possible if the error terms \( \varepsilon_j \) are i.i.d., that is, the noise variances are constant across \( j \)~\citep{peters2014identifiability}. However, the equal variance assumption is not robust to rescaling of the variables.

We conclude this example with a caveat: in the nonlinear case, additive noise models are generally not closed under marginalization, even when marginalizing over unobserved variables. For instance, if \( Y = f_Y(X) + \varepsilon_Y \) and \( Z = f_Z(Y) + \varepsilon_Z \), then the marginal distribution \( P_{X,Z} \) does not necessarily admit an additive noise model from \( X \) to \( Z \). This limitation restricts the applicability of additive noise models, as in most, if not all, scenarios there are intermediate variables on a causal path that are unobserved~\citep[\S 7.1.2]{peters2017elements}. Therefore, additive models should perhaps be viewed as approximations to the true data-generating process rather than exact representations.

% \tblue{In many cases, writing out U can help spelling out identification formula. For example, Ying Zhou's paper; proxy's identification formula. Maybe also Shu Yang's paper? As in this case we can place restrictions on the levels of $U$ through e.g. completeness assumption. cf. pearl chapter 9.4}

% Can also specify the coefficients p.~32

% \tblue{todo: Add a discussion on non-parametric versus parametric (linear) models for causal identification. SEM have parametric version, but should be used with caution. cf. Pearl. p.~94 point 3. cf. pearl chapter 1.4. Can also mention the results in Peters Chapter 4 and 7. And maybe DAGs with no TEARs?}

\section{Summary and Practical Recommendations}

This paper has reviewed three major frameworks for causal inference: the potential outcomes framework, NPSEMs, and DAGs. Each offers a distinct perspective on how to represent, interpret, and identify causal relationships, and each has its own strengths and limitations.

The potential outcomes framework provides a conceptually simple and widely applicable language for defining causal effects, such as average treatment effects and conditional effects. It is especially well suited for settings where study design plays a central role, such as randomized trials and quasi-experiments. Its mathematical simplicity makes it attractive for formulating identification strategies. However, practitioners sometimes  posit counterfactual independence assumptions without providing a substantive justification or circumstances when such assumptions might be likely to hold. These shortcomings can often be addressed using the SWIG/FFRCISTG framework, central to which is the notion of one-step ahead counterfactuals and recursive substitution, which provide a stepwise data-generating process. In the case where there are well-defined interventions on every variable, this leads to a potential outcome model that is isomorphic to an NPSEM.

%the potential outcome lacks a built-in mechanism for encoding complex systems of variables or structural constraints, and assessing the plausibility of identification assumptions.

The NPSEM framework formalizes causal systems using deterministic functional relationships and exogenous error terms. This allows for a generative view of how outcomes arise from their causes, enabling researchers to articulate and evaluate identification assumptions in terms of underlying data-generating mechanisms. As emphasized by \citet[][p.~244]{pearl2009causality}, structural models and their associated graphs are particularly useful as tools for expressing assumptions about cause-effect relationships, as they provide a natural and intuitive representation of mechanistic knowledge. Counterfactual independence assumptions that may appear opaque or difficult to interpret when stated in the potential outcomes notation, such as the ignorability assumption, can often be given immediate and process-based interpretations when expressed structurally. The NPSEM framework also facilitates encoding assumptions on the functional form or even parameter values of the structural functions.  However, when paired with the assumption of independent errors (NPSEM-IE), this framework often imposes strong and untestable cross-world assumptions that go beyond what is needed for identification and remain untestable even under randomized treatment assignment. These added assumptions can constrain the space of admissible models and lead to overly rigid inferences.

The causal DAG framework provides a graphical language for encoding conditional independence assumptions and visualizing the structure of confounding, mediation, and other causal pathways. It facilitates identification via graphical criteria such as the backdoor condition, and its connection to probabilistic DAGs supports a formal link between graphical structure and observed data. Unlike NPSEM-IE, causal DAGs (like SWIGs) avoid  strong cross-world assumptions and are therefore more flexible in some contexts. On the other hand, DAGs do not inherently express counterfactuals, and certain causal quantities, especially those involving nested or cross-world counterfactuals, must be defined outside the graphical formalism or with added assumptions.

Taken together, these frameworks offer complementary tools for causal inference. A practical and robust workflow often benefits from combining them:
\begin{enumerate}
    \item Define the causal estimand in potential outcomes notation to clarify the scientific question. This step requires basic assumptions such as the SUTVA, but does not otherwise rely on parametric modeling assumptions.
    \item Express and evaluate assumptions using DAGs, SWIGs or NPSEMs to articulate structural knowledge and identify valid adjustment strategies.
    \item Establish identification using graphical criteria such as d-separation, 
and algebraic tools like the g-formula, do-calculus, po-calculus or ID algorithm.
    \item Choose an estimation method appropriate to the structure and available data, drawing on statistical tools such as regression, weighting, or machine learning methods.
\end{enumerate}

% This integrated approach not only clarifies which assumptions are essential for causal identification, but also helps avoid imposing unnecessary or implausible constraints. Structural models and DAGs offer transparency and intuition, while potential outcomes provide a precise language for formulating causal questions and guiding estimation.

In this workflow, structural assumptions are naturally expressed in graphical or NPSEM terms, then translated into potential outcomes notation, and finally subjected to algebraic derivation \citep[][p.~245]{pearl2009causality} to isolate those that are strictly necessary for identification. In this way, NPSEMs and their associated graphs (see Section~\ref{sec:npsem-dag} for details) serve as tools for evaluating the plausibility or implausibility of conditional independence assumptions stated in terms of potential outcomes. Contemporary implementations such as the target trial emulation framework \citep{hernan2016using, hernan2022target} illustrate how structured causal roadmaps can help prevent common inferential pitfalls in observational studies. We refer readers to \citet{wang2017identification,wang2018bounded,yang2019causal,zhou2024promises,dong2025marginal} for examples that follow this workflow.

We now illustrate these steps in detail using \citet{wang2018bounded}. Their study begins by defining the causal estimand, the ACE, in potential outcomes notation, thereby clarifying the scientific target of inference (Step 1). They then express structural assumptions through a DAG and a corresponding SWIG with bidirected arrows, explicitly distinguishing identification assumptions such as the exclusion restriction and independence from modeling assumptions (Step 2). Using algebraic derivations, they establish identification of the ACE under new no-interaction assumptions (A5.a or A5.b), showing that the estimand can be represented as the average Wald estimand (Step 3). Finally, they develop a series of bounded, efficient, and multiply robust estimators that remain consistent if any one of several working models is correctly specified (Step 4). This example demonstrates how potential outcomes, graphical models, and semiparametric theory can be systematically combined into a coherent causal workflow.

Recent advances continue to expand the reach of causal inference beyond traditional domains. Active research areas include causal discovery algorithms for learning causal structure from data \citep[e.g.,][]{spirtes2000causation,peters2017elements}, fairness in artificial intelligence \citep[e.g.,][]{kusner2017counterfactual}, and methods for robust learning and generalization across heterogeneous environments \citep[e.g.,][]{peters2016causal}. Related ideas have also influenced reinforcement learning and decision-making in artificial intelligence, where causal reasoning is increasingly recognized as essential for interpretability and transferability. A full treatment of these topics is beyond the scope of this review, but they exemplify how causal thinking continues to shape contemporary research at the intersection of statistics, computer science, and the social and health sciences.

\section*{Acknowledgements}
% The conceptual foundation of this paper originates from a lecture delivered by Thomas Richardson during the author's graduate studies. 
We are grateful to the Editor for encouraging the development of this work, and to Forrest Crawford  and Sander Greenland for insightful discussions that helped shape several key aspects. 
% We are especially grateful to Thomas Richardson for carefully reading the draft and providing constructive feedback. 
We acknowledge the use of OpenAI's ChatGPT for proofreading and refining the writing.

 \section*{Appendix}

\section*{A1. Constructing Priors and “Sampling Uniformly’’ on $\mathcal S_2$}

We now provide a principled way to “sample uniformly’’ from the space
\[
\mathcal S_2 = \{ P : A \ind B,\ A \ind D \}.
\]
The strategy is to express $\mathcal S_2$ using a set of variation independent
parameters. This creates a simple product space on which a prior can be placed, and
the resulting measure pushes forward to a well-defined prior on $\mathcal S_2$.

Throughout we write
\[
p_a^A = P(A=a),\quad
p_b^B = P(B=b),\quad
p_d^D = P(D=d),
\qquad a,b,d\in\{0,1\}.
\]
The assumptions $A\ind B$ and $A\ind D$ imply that the conditional marginals of
$B$ and $D$ do not depend on $A$, so
\[
P(B=b\mid A=a)=p_b^B,\qquad P(D=d\mid A=a)=p_d^D
\]
for all $a$. What remains is to characterize the association between $B$ and $D$
within each stratum of $A$, while preserving these common one way marginals.

For binary variables, it is well known that a $2\times2$ table with fixed marginals
has exactly one degree of freedom. Any strictly positive such table can therefore be
uniquely represented by a single scalar measuring the strength of association.
A standard choice for this scalar is the conditional odds ratio
\citep{chen2007semiparametric,tchetgen2018discrete}:
\[
\mathrm{OR}(B,D \mid A=a)
  = \frac{P(B=1,D=1\mid A=a)\,P(B=0,D=0\mid A=a)}
         {P(B=1,D=0\mid A=a)\,P(B=0,D=1\mid A=a)}.
\]
The odds ratio is variation independent of the one way marginals of $B$ and $D$, and
under $A\ind B$ and $A\ind D$, these conditional marginals coincide with the
unconditional marginals. Thus, specifying $\mathrm{OR}(B,D\mid A=a)$ together with
$(p_1^B,p_1^D)$ fully determines the entire conditional $2\times2$ table in stratum
$A=a$.

These observations motivate parameterizing $\mathcal S_2$ by the three free marginal
probabilities $(p_1^A,p_1^B,p_1^D)\in(0,1)^3$ together with the two stratum specific
odds ratios
\[
\psi_a=\mathrm{OR}(B,D\mid A=a)\in(0,\infty),\qquad a=0,1.
\]
The five parameters $(p_1^A,p_1^B,p_1^D,\psi_0,\psi_1)$ vary freely over the Cartesian
product $(0,1)^3\times(0,\infty)^2$ and uniquely determine a strictly positive
distribution in $\mathcal S_2$.

To construct a prior on $\mathcal S_2$, one may place independent
$\mathrm{Dirichlet}(1,1)$ priors on $(p_0^A,p_1^A)$, $(p_0^B,p_1^B)$ and
$(p_0^D,p_1^D)$, together with independent priors on the log odds ratios, for example
\[
\log \psi_0,\ \log \psi_1 \sim \mathrm{Normal}(0,1).
\]
 This induces a smooth prior with full support on the strictly
positive interior of $\mathcal S_2$ and provides a concrete interpretation of
“sampling uniformly’’ over $\mathcal S_2$.

\bibliography{causal}

\begin{thebibliography}{}

\bibitem[Aldrich, 1989]{aldrich1989autonomy}
Aldrich, J. (1989).
\newblock Autonomy.
\newblock {\em Oxford Economic Papers}, 41(1):15--34.

\bibitem[Angrist et~al., 1996]{angrist1996identification}
Angrist, J.~D., Imbens, G.~W., and Rubin, D.~B. (1996).
\newblock Identification of causal effects using instrumental variables.
\newblock {\em Journal of the American Statistical Association}, 91(434):444--455.

\bibitem[Arjovsky et~al., 2020]{Arjovsky2020IRM}
Arjovsky, M., Bottou, L., Gulrajani, I., and Lopez-Paz, D. (2020).
\newblock Invariant risk minimization.
\newblock {\em arXiv preprint arXiv:1907.02893}.
\newblock Original version 2019.

\bibitem[Bajari et~al., 2023]{bajari2023experimental}
Bajari, P., Burdick, B., Imbens, G.~W., Masoero, L., McQueen, J., Richardson, T.~S., and Rosen, I.~M. (2023).
\newblock Experimental design in marketplaces.
\newblock {\em Statistical Science}, 38(3):458--476.

\bibitem[Balke and Pearl, 1997]{balke1997nonparametric}
Balke, A. and Pearl, J. (1997).
\newblock Bounds on treatment effects from studies with imperfect compliance.
\newblock {\em Journal of the American Statistical Association}, 92(439):1171--1176.

\bibitem[Billingsley, 1995]{billingsley1995probability}
Billingsley, P. (1995).
\newblock {\em Probability and Measure}.
\newblock Wiley, New York, 3rd edition.

\bibitem[B{\"u}hlmann, 2020]{Buhlmann2020StatSci}
B{\"u}hlmann, P. (2020).
\newblock Invariance, causality and robustness.
\newblock {\em Statistical Science}, 35(3):404--426.

\bibitem[Cartwright, 2007]{cartwright2007hunting}
Cartwright, N. (2007).
\newblock {\em Hunting Causes and Using Them: Approaches in Philosophy and Economics}.
\newblock Cambridge University Press, Cambridge, UK.

\bibitem[Chen, 2007]{chen2007semiparametric}
Chen, H.~Y. (2007).
\newblock A semiparametric odds ratio model for measuring association.
\newblock {\em Biometrics}, 63(2):413--421.

\bibitem[Cole and Frangakis, 2009]{cole2009consistency}
Cole, S.~R. and Frangakis, C.~E. (2009).
\newblock The consistency statement in causal inference: {A} definition or an assumption?
\newblock {\em Epidemiology}, 20(1):3--5.

\bibitem[Dawid, 2000]{dawid2000causal}
Dawid, A.~P. (2000).
\newblock Causal inference without counterfactuals.
\newblock {\em Journal of the American Statistical Association}, 95(450):407--424.

\bibitem[Dawid, 2015]{dawid2015statistical}
Dawid, A.~P. (2015).
\newblock Statistical causality from a decision-theoretic perspective.
\newblock {\em Annual Review of Statistics and Its Application}, 2(1):273--303.

\bibitem[Ding et~al., 2011]{ding2011identifiability}
Ding, P., Geng, Z., Yan, W., and Zhou, X.-H. (2011).
\newblock Identifiability and estimation of causal effects by principal stratification with outcomes truncated by death.
\newblock {\em Journal of the American Statistical Association}, 106(496):1578--1591.

\bibitem[Dong et~al., 2025]{dong2025marginal}
Dong, M., Liu, L., Tang, D., Liu, G., Xu, W., and Wang, L. (2025).
\newblock Marginal causal effect estimation with continuous instrumental variables.
\newblock {\em arXiv preprint arXiv:2510.14368}.

\bibitem[Druzdzel and Simon, 1993]{druzdzel1993causality}
Druzdzel, M.~J. and Simon, H.~A. (1993).
\newblock Causality in {B}ayesian belief networks.
\newblock In {\em Proceedings of the Ninth Conference on Uncertainty in Artificial Intelligence (UAI-93)}, pages 3--11, San Francisco, CA. Morgan Kaufmann.

\bibitem[El~Gamal and Kim, 2011]{el2011network}
El~Gamal, A. and Kim, Y.-H. (2011).
\newblock {\em Network {I}nformation {T}heory}.
\newblock Cambridge University Press.

\bibitem[Geiger et~al., 1990]{geiger1990identifying}
Geiger, D., Verma, T., and Pearl, J. (1990).
\newblock Identifying independence in {B}ayesian networks.
\newblock {\em Networks}, 20(5):507--534.

\bibitem[Greenland, 2004]{Greenland2004Overview}
Greenland, S. (2004).
\newblock An overview of methods for causal inference from observational studies.
\newblock In Gelman, A. and Meng, X.-L., editors, {\em Applied Bayesian Modeling and Causal Inference from Incomplete-Data Perspectives: An Essential Journey with Donald Rubin's Statistical Family}, pages 1--13. John Wiley \& Sons, Chichester, UK.

\bibitem[Greenland and Brumback, 2002]{greenland2002overview}
Greenland, S. and Brumback, B. (2002).
\newblock An overview of relations among causal modelling methods.
\newblock {\em International Journal of Epidemiology}, 31(5):1030--1037.

\bibitem[Greenland and Poole, 1988]{greenland1988invariants}
Greenland, S. and Poole, C. (1988).
\newblock Invariants and noninvariants in the concept of interdependent effects.
\newblock {\em Scandinavian Journal of Work, Environment \& Health}, 14(2):125--129.

\bibitem[Haavelmo, 1943]{haavelmo1943probability}
Haavelmo, T. (1943).
\newblock The statistical implications of a system of simultaneous equations.
\newblock {\em Econometrica}, 11(1):1--12.

\bibitem[Halpern and Pearl, 2000]{halpern2000axiomatizing}
Halpern, J.~Y. and Pearl, J. (2000).
\newblock Axiomatizing causal reasoning.
\newblock {\em Journal of Artificial Intelligence Research}, 12:317--337.

\bibitem[Hartwig et~al., 2023]{hartwig2023average}
Hartwig, F.~P., Wang, L., Smith, G.~D., and Davies, N.~M. (2023).
\newblock Average causal effect estimation via instrumental variables: {T}he no simultaneous heterogeneity assumption.
\newblock {\em Epidemiology}, 34(3):325--332.

\bibitem[Hern{\'a}n and Robins, 2016]{hernan2016using}
Hern{\'a}n, M.~A. and Robins, J.~M. (2016).
\newblock Using big data to emulate a target trial when a randomized trial is not available.
\newblock {\em American Journal of Epidemiology}, 183(8):758--764.

\bibitem[Hern{\'a}n and Robins, 2025]{hernan2020causal}
Hern{\'a}n, M.~A. and Robins, J.~M. (2025).
\newblock {\em Causal Inference: What If}.
\newblock Chapman \& Hall/CRC, Boca Raton, FL, 1st edition.

\bibitem[Hern{\'a}n et~al., 2022]{hernan2022target}
Hern{\'a}n, M.~A., Wang, W., and Leaf, D.~E. (2022).
\newblock Target trial emulation: a framework for causal inference from observational data.
\newblock {\em Journal of the American Medical Association}, 328(24):2446--2447.

\bibitem[Hill, 1965]{Hill1965Environment}
Hill, A.~B. (1965).
\newblock The environment and disease: Association or causation?
\newblock {\em Proceedings of the Royal Society of Medicine}, 58:295--300.

\bibitem[Holland, 1986]{holland1986statistics}
Holland, P.~W. (1986).
\newblock Statistics and causal inference.
\newblock {\em Journal of the American Statistical Association}, 81(396):945--960.

\bibitem[Hoyer et~al., 2008]{hoyer2008nonlinear}
Hoyer, P., Janzing, D., Mooij, J.~M., Peters, J., and Sch\"{o}lkopf, B. (2008).
\newblock Nonlinear causal discovery with additive noise models.
\newblock In Koller, D., Schuurmans, D., Bengio, Y., and Bottou, L., editors, {\em Advances in Neural Information Processing Systems}, volume~21. Curran Associates, Inc.

\bibitem[Huang and Valtorta, 2006]{huang2006pearl}
Huang, Y. and Valtorta, M. (2006).
\newblock Pearl's calculus of intervention is complete.
\newblock In {\em Proceedings of the 22nd Conference on Uncertainty in Artificial Intelligence (UAI-2006)}, pages 217--224. AUAI Press.

\bibitem[Imai et~al., 2010]{imai2010identification}
Imai, K., Keele, L., and Yamamoto, T. (2010).
\newblock Identification, inference and sensitivity analysis for causal mediation effects.
\newblock {\em Statistical Science}, 25(1):51--71.

\bibitem[Imbens, 2014]{imbens:2014}
Imbens, G.~W. (2014).
\newblock Instrumental variables: An econometrician’s perspective.
\newblock {\em Statistical Science}, 29(3):323 -- 358.

\bibitem[Imbens and Rubin, 2015]{imbens2015causal}
Imbens, G.~W. and Rubin, D.~B. (2015).
\newblock {\em Causal Inference for Statistics, Social, and Biomedical Sciences: An Introduction}.
\newblock Cambridge University Press, New York.

\bibitem[Imbens and Wooldridge, 2009]{ImbensWooldridge2009JEL}
Imbens, G.~W. and Wooldridge, J.~M. (2009).
\newblock Recent developments in the econometrics of program evaluation.
\newblock {\em Journal of Economic Literature}, 47(1):5--86.

\bibitem[Janzing and Schl{\"o}lkopf, 2010]{janzing2010causal}
Janzing, D. and Schl{\"o}lkopf, B. (2010).
\newblock Causal inference using the algorithmic {Markov} condition.
\newblock {\em IEEE Transactions on Information Theory}, 56(10):5168--5194.

\bibitem[Jiao et~al., 2024]{Jiao2024ResearchSurvey}
Jiao, L., Wang, Y., Liu, X., Li, L., Liu, F., Ma, W., Guo, Y., Chen, P., Yang, S., and Hou, B. (2024).
\newblock Causal inference meets deep learning: A comprehensive survey.
\newblock {\em Research}, 7:0467.

\bibitem[Jun and Lee, 2023]{jun2023identifying}
Jun, S.~J. and Lee, S. (2023).
\newblock Identifying the effect of persuasion.
\newblock {\em Journal of Political Economy}, 131(8):2032--2058.

\bibitem[Kuang et~al., 2020]{kuang2020causal}
Kuang, K., Li, L., Geng, Z., Xu, L., Zhang, K., Liao, B., Huang, H., Ding, P., Miao, W., and Jiang, Z. (2020).
\newblock Causal inference.
\newblock {\em Engineering}, 6:253--263.

\bibitem[Kusner et~al., 2017]{kusner2017counterfactual}
Kusner, M.~J., Loftus, J., Russell, C., and Silva, R. (2017).
\newblock Counterfactual fairness.
\newblock {\em Advances in Neural Information Processing Systems}, 30.

\bibitem[Laplace, 1825]{laplace1825essai}
Laplace, P. (1825).
\newblock {\em Essai Philosophique sur les Probabilit{\'e}s}.
\newblock Bachelier, Paris, 5 edition.
\newblock Originally published in 1814.

\bibitem[Lauritzen, 1996]{lauritzen1996graphical}
Lauritzen, S.~L. (1996).
\newblock {\em Graphical Models}, volume~17.
\newblock Clarendon Press.

\bibitem[Lauritzen, 2004]{lauritzen2004discussion}
Lauritzen, S.~L. (2004).
\newblock Discussion on causality. {D}iscussion of {A. P. D}awid's ``{Probability, causality, and the empirical world}: A {B}ayes--de {F}inetti--{P}opper--{B}orel synthesis''.
\newblock {\em Scandinavian Journal of Statistics}, 31(2):189--192.

\bibitem[Lauritzen et~al., 1990]{lauritzen1990independence}
Lauritzen, S.~L., Dawid, A.~P., Larsen, B.~N., and Leimer, H.-G. (1990).
\newblock Independence properties of directed {Markov} fields.
\newblock {\em Networks}, 20(5):491--505.

\bibitem[Lehmann and Casella, 1998]{lehmann2006theory}
Lehmann, E.~L. and Casella, G. (1998).
\newblock {\em Theory of Point Estimation}.
\newblock Springer, New York, 2 edition.

\bibitem[Li et~al., 2023]{li2023bayesian}
Li, F., Ding, P., and Mealli, F. (2023).
\newblock Bayesian causal inference: A critical review.
\newblock {\em Philosophical Transactions of the Royal Society A}, 381(2247):20220153.

\bibitem[Malinsky et~al., 2019]{malinsky-19b}
Malinsky, D., Shpitser, I., and Richardson, T. (2019).
\newblock A potential outcomes calculus for identifying conditional path-specific effects.
\newblock In Chaudhuri, K. and Sugiyama, M., editors, {\em Proceedings of the Twenty-Second International Conference on Artificial Intelligence and Statistics}, volume~89 of {\em Proceedings of Machine Learning Research}, pages 3080--3088. PMLR.

\bibitem[Manski, 1993]{manski1999identification}
Manski, C.~F. (1993).
\newblock Identification problems in the social sciences.
\newblock {\em Sociological Methodology}, 23:1--56.

\bibitem[Morgan and Winship, 2014]{morgan2014counterfactuals}
Morgan, S.~L. and Winship, C. (2014).
\newblock {\em Counterfactuals and Causal Inference: Methods and Principles for Social Research}.
\newblock Cambridge University Press, Cambridge, 2 edition.

\bibitem[Neison, 1844]{neison1844method}
Neison, F. (1844).
\newblock On a method recently proposed for conducting inquiries into the comparative sanatory condition of various districts, with illustrations, derived from numerous places in {G}reat {B}ritain at the period of the last census.
\newblock {\em Journal of the Statistical Society of London}, 7(1):40--68.

\bibitem[Neyman, 1923]{neyman1923application}
Neyman, J. (1923).
\newblock On the application of probability theory to agricultural experiments. {E}ssay on principles. {S}ection 9.
\newblock {\em Statistical Science}, 5(4):465--472.
\newblock Translated and edited by D. M. Dabrowska and T. P. Speed from the 1923 Polish original.

\bibitem[Pearl, 1985]{pearl1985bayesian}
Pearl, J. (1985).
\newblock Bayesian networks: A model of self-activated memory for evidential reasoning.
\newblock {\em Proceedings of the 7th Conference of the Cognitive Science Society}, pages 329--334.

\bibitem[Pearl, 1988]{pearl1988probabilistic}
Pearl, J. (1988).
\newblock {\em Probabilistic Reasoning in Intelligent Systems: Networks of Plausible Inference}.
\newblock Morgan Kaufmann, San Mateo, CA.

\bibitem[Pearl, 1995]{pearl1995causal}
Pearl, J. (1995).
\newblock Causal diagrams for empirical research.
\newblock {\em Biometrika}, 82(4):669--688.

\bibitem[Pearl, 2001]{pearl2001direct}
Pearl, J. (2001).
\newblock Direct and indirect effects.
\newblock In Breese, J.~S. and Koller, D., editors, {\em Proceedings of the Seventeenth Conference on Uncertainty in Artificial Intelligence (UAI-2001)}, pages 411--420, San Francisco, CA. Morgan Kaufmann.

\bibitem[Pearl, 2009]{pearl2009causality}
Pearl, J. (2009).
\newblock {\em Causality}.
\newblock Cambridge University Press, Cambridge, 2nd ed., first printing edition.

\bibitem[Pearl, 2010a]{pearl2010IJB}
Pearl, J. (2010a).
\newblock An introduction to causal inference.
\newblock {\em The International Journal of Biostatistics}, 6(2):7.

\bibitem[Pearl, 2010b]{pearl2010consistency}
Pearl, J. (2010b).
\newblock On the consistency rule in causal inference: Axiom, definition, assumption, or theorem?
\newblock {\em Epidemiology}, 21(6):872--875.

\bibitem[Pearl, 2019]{Pearl2019CACMSevenTools}
Pearl, J. (2019).
\newblock The seven tools of causal inference, with reflections on machine learning.
\newblock {\em Communications of the ACM}, 62(3):54--60.

\bibitem[Pearl et~al., 2016]{pearl2016causal}
Pearl, J., Glymour, M., and Jewell, N.~P. (2016).
\newblock {\em Causal Inference in Statistics: A Primer}.
\newblock John Wiley \& Sons.

\bibitem[Pearl and Mackenzie, 2018]{pearl2018bookofwhy}
Pearl, J. and Mackenzie, D. (2018).
\newblock {\em The Book of Why: The New Science of Cause and Effect}.
\newblock Basic Books, New York.

\bibitem[Pearl and Paz, 1986]{pearl1986graphoids}
Pearl, J. and Paz, A. (1986).
\newblock Graphoids: A graph-based logic for reasoning about relevance relations.
\newblock In Boulay, B.~D., Warren, D. H.~D., and Kyburg, H.~E., editors, {\em Advances in Artificial Intelligence II}, pages 357--363. North-Holland.

\bibitem[Perkovi{\'c} et~al., 2018]{Perkovic2018CompleteAdjustment}
Perkovi{\'c}, E., Textor, J., Kalisch, M., and Maathuis, M.~H. (2018).
\newblock Complete graphical characterization and construction of adjustment sets in markov equivalence classes of ancestral graphs.
\newblock {\em Journal of Machine Learning Research}, 18(220):1--62.

\bibitem[Peters and B{\"u}hlmann, 2014]{peters2014identifiability}
Peters, J. and B{\"u}hlmann, P. (2014).
\newblock Identifiability of {G}aussian structural equation models with equal error variances.
\newblock {\em Biometrika}, 101(1):219--228.

\bibitem[Peters et~al., 2016]{peters2016causal}
Peters, J., B{\"u}hlmann, P., and Meinshausen, N. (2016).
\newblock Causal inference by using invariant prediction: Identification and confidence intervals.
\newblock {\em Journal of the Royal Statistical Society Series B: Statistical Methodology}, 78(5):947--1012.

\bibitem[Peters et~al., 2017]{peters2017elements}
Peters, J., Janzing, D., and Sch{\"o}lkopf, B. (2017).
\newblock {\em Elements of Causal Inference: Foundations and Learning Algorithms}.
\newblock The MIT Press.

\bibitem[Peters et~al., 2014]{peters2014causal}
Peters, J., Mooij, J.~M., Janzing, D., and Sch{\"o}lkopf, B. (2014).
\newblock Causal discovery with continuous additive noise models.
\newblock {\em Journal of Machine Learning Research}, 15(1):2009--2053.

\bibitem[Richardson, 2003]{richardson2003markov}
Richardson, T.~S. (2003).
\newblock Markov properties for acyclic directed mixed graphs.
\newblock {\em Scandinavian Journal of Statistics}, 30(1):145--157.

\bibitem[Richardson et~al., 2023]{richardson2023nested}
Richardson, T.~S., Evans, R.~J., Robins, J.~M., and Shpitser, I. (2023).
\newblock Nested {M}arkov properties for acyclic directed mixed graphs.
\newblock {\em The Annals of Statistics}, 51(1):334--361.

\bibitem[Richardson and Robins, 2013]{richardson2013single}
Richardson, T.~S. and Robins, J.~M. (2013).
\newblock Single world intervention graphs ({SWIGs}): A unification of the counterfactual and graphical approaches to causality.
\newblock Technical Report Working Paper 128, Center for Statistics and the Social Sciences, University of Washington.

\bibitem[Richardson and Robins, 2014]{richardson:robins:2014}
Richardson, T.~S. and Robins, J.~M. (2014).
\newblock {ACE} bounds; {SEMs} with equilibrium conditions.
\newblock {\em Statistical Science}, 29(3):363--366.

\bibitem[Richardson and Robins, 2023]{richardson2023potential}
Richardson, T.~S. and Robins, J.~M. (2023).
\newblock Potential outcome and decision theoretic foundations for statistical causality.
\newblock {\em Journal of Causal Inference}, 11(1):20220012.

\bibitem[Richardson and Spirtes, 2002]{richardson2002ancestral}
Richardson, T.~S. and Spirtes, P. (2002).
\newblock Ancestral graph {Markov} models.
\newblock {\em The Annals of Statistics}, 30(4):962--1030.

\bibitem[Robins and Greenland, 1989]{robins1989probability}
Robins, J. and Greenland, S. (1989).
\newblock The probability of causation under a stochastic model for individual risk.
\newblock {\em Biometrics}, 45(4):1125--1138.

\bibitem[Robins and Greenland, 1991]{robins1991estimability}
Robins, J. and Greenland, S. (1991).
\newblock Estimability and estimation of expected years of life lost due to a hazardous exposure.
\newblock {\em Statistics in Medicine}, 10(1):79--93.

\bibitem[Robins, 1986]{robins1986new}
Robins, J.~M. (1986).
\newblock A new approach to causal inference in mortality studies with a sustained exposure period—application to control of the healthy worker survivor effect.
\newblock {\em Mathematical Modelling}, 7(9–12):1393--1512.

\bibitem[Robins, 2003]{robins2003semantics}
Robins, J.~M. (2003).
\newblock Semantics of causal {DAG} models and the identification of direct and indirect effects.
\newblock In Green, P.~J., Hjort, N.~L., and Richardson, S., editors, {\em Highly Structured Stochastic Systems}, pages 70--81. Oxford University Press, Oxford.

\bibitem[Robins and Greenland, 1992]{robins:greenland:1992}
Robins, J.~M. and Greenland, S. (1992).
\newblock Identifiability and exchangeability for direct and indirect effects.
\newblock {\em Epidemiology}, 3(2):143--155.

\bibitem[Robins et~al., 2004]{robins2004effects}
Robins, J.~M., Hern{\'a}n, M.~A., and Siebert, U. (2004).
\newblock Effects of multiple interventions.
\newblock In Ezzati, M., Lopez, A.~D., Rodgers, A., and Murray, C. J.~L., editors, {\em Comparative Quantification of Health Risks: Global and Regional Burden of Disease Attributable to Selected Major Risk Factors}, pages 2191--2230. World Health Organization, Geneva.

\bibitem[Robins and Richardson, 2011]{robins2011alternative}
Robins, J.~M. and Richardson, T.~S. (2011).
\newblock Alternative graphical causal models and the identification of direct effects.
\newblock In Shrout, P., Keyes, K., and Ornstein, K., editors, {\em Causality and Psychopathology: Finding the Determinants of Disorders and their Cures}, chapter~6, pages 1--52. Oxford University Press.

\bibitem[Rosenbaum, 2010]{rosenbaum2010design}
Rosenbaum, P.~R. (2010).
\newblock {\em Design of Observational Studies}.
\newblock Springer Series in Statistics. Springer, New York.

\bibitem[Rothman, 1976]{rothman1976causes}
Rothman, K.~J. (1976).
\newblock Causes.
\newblock {\em American Journal of Epidemiology}, 104(6):587--592.

\bibitem[Rothman et~al., 2008]{rothman2008modern}
Rothman, K.~J., Greenland, S., Lash, T.~L., et~al. (2008).
\newblock {\em Modern Epidemiology}, volume~3.
\newblock Wolters Kluwer Health/Lippincott Williams \& Wilkins Philadelphia.

\bibitem[Roy, 1951]{roy1951some}
Roy, A.~D. (1951).
\newblock Some thoughts on the distribution of earnings.
\newblock {\em Oxford Economic Papers}, 3(2):135--146.

\bibitem[Rubin, 1974]{rubin1974estimating}
Rubin, D.~B. (1974).
\newblock Estimating causal effects of treatments in randomized and nonrandomized studies.
\newblock {\em Journal of Educational Psychology}, 66(5):688--701.

\bibitem[Rubin, 1980]{rubin1980randomization}
Rubin, D.~B. (1980).
\newblock Randomization analysis of experimental data: The {F}isher randomization test comment.
\newblock {\em Journal of the American Statistical Association}, 75(371):591--593.

\bibitem[Rubin, 2004]{rubin2004direct}
Rubin, D.~B. (2004).
\newblock Direct and indirect causal effects via potential outcomes.
\newblock {\em Scandinavian Journal of Statistics}, 31(2):161--170.

\bibitem[Rubin, 2005]{rubin2005causal}
Rubin, D.~B. (2005).
\newblock Causal inference using potential outcomes: Design, modeling, decisions.
\newblock {\em Journal of the American Statistical Association}, 100(469):322--331.

\bibitem[Sch{\"o}lkopf, 2022]{Scholkopf2022CausalityML}
Sch{\"o}lkopf, B. (2022).
\newblock Causality for machine learning.
\newblock In Geffner, H., Dechter, R., and Halpern, J.~Y., editors, {\em Probabilistic and Causal Inference: The Works of Judea Pearl}, ACM Books, pages 765--804. Association for Computing Machinery, New York, NY.

\bibitem[Shimizu et~al., 2006]{shimizu2006linear}
Shimizu, S., Hoyer, P.~O., Hyv{\"{a}}rinen, A., and Kerminen, A. (2006).
\newblock A linear non-{G}aussian acyclic model for causal discovery.
\newblock {\em Journal of Machine Learning Research}, 7(72):2003--2030.

\bibitem[Shpitser and Pearl, 2006]{shpitser2006identification}
Shpitser, I. and Pearl, J. (2006).
\newblock Identification of joint interventional distributions in recursive {semi-Markovian} causal models.
\newblock In {\em AAAI}, pages 1219--1226.

\bibitem[Shpitser and Pearl, 2007]{shpitser2007what}
Shpitser, I. and Pearl, J. (2007).
\newblock What counterfactuals can be tested.
\newblock In {\em Proceedings of the 23rd Conference on Uncertainty in Artificial Intelligence (UAI)}, pages 352--359. AUAI Press.

\bibitem[Shpitser and Pearl, 2008]{shpitser2008complete}
Shpitser, I. and Pearl, J. (2008).
\newblock Complete identification methods for the causal hierarchy.
\newblock {\em Journal of Machine Learning Research}, 9:1941--1979.

\bibitem[Shpitser et~al., 2022]{shpitser2022multivariate}
Shpitser, I., Richardson, T.~S., and Robins, J.~M. (2022).
\newblock Multivariate counterfactual systems and causal graphical models.
\newblock In {\em Probabilistic and Causal Inference: The Works of Judea Pearl}, pages 813--852.

\bibitem[Simon, 1953]{simon1953causal}
Simon, H.~A. (1953).
\newblock {\em Causal ordering and identifiability}.
\newblock Studies in Econometric Method. Wiley, New York.

\bibitem[Spirtes, 2010]{spirtes2010introduction}
Spirtes, P. (2010).
\newblock Introduction to causal inference.
\newblock {\em Journal of Machine Learning Research}, 11:1643--1662.

\bibitem[Spirtes et~al., 2000]{spirtes2000causation}
Spirtes, P., Glymour, C.~N., and Scheines, R. (2000).
\newblock {\em Causation, Prediction, and Search}.
\newblock MIT Press, Cambridge, MA, 2nd edition.

\bibitem[Strotz and Wold, 1960]{strotz1960recursive}
Strotz, R.~H. and Wold, H. O.~A. (1960).
\newblock Recursive and nonrecursive systems of equations.
\newblock {\em Econometrica}, 28(2):417--427.

\bibitem[Tchetgen~Tchetgen et~al., 2018]{tchetgen2018discrete}
Tchetgen~Tchetgen, E.~J., Wang, L., and Sun, B. (2018).
\newblock Discrete choice models for nonmonotone nonignorable missing data: Identification and inference.
\newblock {\em Statistica Sinica}, 28(4):2069.

\bibitem[Tian and Pearl, 2002]{tian2002identification}
Tian, J. and Pearl, J. (2002).
\newblock A general identification condition for causal effects.
\newblock In {\em Proceedings of the Eighteenth National Conference on Artificial Intelligence (AAAI-2002)}, pages 567--573. AAAI Press.

\bibitem[Tian and Shpitser, 2010]{tian2010identifying}
Tian, J. and Shpitser, I. (2010).
\newblock On identifying causal effects.
\newblock {\em Heuristics, Probability and Causality: A Tribute to Judea Pearl (R. Dechter, H. Geffner and J. Halpern, eds.). College Publications, UK}, pages 415--444.

\bibitem[Tjoa and Guan, 2021]{TjoaGuan2021TNNLSXAIHealth}
Tjoa, E. and Guan, C. (2021).
\newblock A survey on explainable artificial intelligence ({XAI}): Toward medical {XAI}.
\newblock {\em IEEE Transactions on Neural Networks and Learning Systems}, 32(11):4793--4813.

\bibitem[Uhler et~al., 2013]{uhler2013geometry}
Uhler, C., Raskutti, G., B{\"u}hlmann, P., and Yu, B. (2013).
\newblock Geometry of the faithfulness assumption in causal inference.
\newblock {\em The Annals of Statistics}, 41(2):436--463.

\bibitem[VanderWeele, 2015]{vanderweele2015explanation}
VanderWeele, T. (2015).
\newblock {\em Explanation in Causal Inference: Methods for Mediation and Interaction}.
\newblock Oxford University Press.

\bibitem[VanderWeele, 2009]{vanderweele2009concerning}
VanderWeele, T.~J. (2009).
\newblock Concerning the consistency assumption in causal inference.
\newblock {\em Epidemiology}, 20(6):880--883.

\bibitem[VanderWeele and Richardson, 2012]{vanderweele2012general}
VanderWeele, T.~J. and Richardson, T.~S. (2012).
\newblock General theory for interactions in sufficient cause models with dichotomous exposures.
\newblock {\em The Annals of Statistics}, 40(4):2128.

\bibitem[VanderWeele and Robins, 2009]{vanderweele2009minimal}
VanderWeele, T.~J. and Robins, J.~M. (2009).
\newblock Minimal sufficient causation and directed acyclic graphs.
\newblock {\em The Annals of Statistics}, 37(3):1437--1465.

\bibitem[Verma and Pearl, 1991]{verma1991equivalence}
Verma, T.~S. and Pearl, J. (1991).
\newblock Equivalence and synthesis of causal models.
\newblock In {\em Proceedings of the Sixth Conference on Uncertainty in Artificial Intelligence (UAI)}, pages 220--227. Morgan Kaufmann.

\bibitem[Wachter et~al., 2017]{Wachter2017HJLTCounterfactual}
Wachter, S., Mittelstadt, B., and Russell, C. (2017).
\newblock Counterfactual explanations without opening the black box: Automated decisions and the {GDPR}.
\newblock {\em Harvard Journal of Law \& Technology}, 31(2):841--887.

\bibitem[Wang, 2022]{wang2022homogeneity}
Wang, L. (2022).
\newblock On the homogeneity of measures for binary associations.
\newblock {\em arXiv preprint arXiv:2210.05179}.

\bibitem[Wang et~al., 2017a]{wang2017causal}
Wang, L., Richardson, T.~S., and Zhou, X.-H. (2017a).
\newblock Causal analysis of ordinal treatments and binary outcomes under truncation by death.
\newblock {\em Journal of the Royal Statistical Society: Series B (Statistical Methodology)}, 79(3):719--735.

\bibitem[Wang and Tchetgen~Tchetgen, 2018]{wang2018bounded}
Wang, L. and Tchetgen~Tchetgen, E.~J. (2018).
\newblock Bounded, efficient and multiply robust estimation of average treatment effects using instrumental variables.
\newblock {\em Journal of the Royal Statistical Society: Series B (Statistical Methodology)}, 80(3):531--550.

\bibitem[Wang et~al., 2017b]{wang2017identification}
Wang, L., Zhou, X.-H., and Richardson, T.~S. (2017b).
\newblock Identification and estimation of causal effects with outcomes truncated by death.
\newblock {\em Biometrika}, 104(3):597--612.

\bibitem[Wright, 1921]{wright1921correlation}
Wright, S. (1921).
\newblock Correlation and causation.
\newblock {\em Journal of Agricultural Research}, 20:557--585.

\bibitem[Wu and Wang, 2026]{wu2026position}
Wu, P. and Wang, L. (2026).
\newblock Position: A potential outcomes perspective on {P}earl's causal hierarchy.
\newblock {\em arXiv preprint arXiv:2601.20405}.

\bibitem[Yang et~al., 2019]{yang2019causal}
Yang, S., Wang, L., and Ding, P. (2019).
\newblock Causal inference with confounders missing not at random.
\newblock {\em Biometrika}, 106(4):875--888.

\bibitem[Yao et~al., 2021]{Yao2021Survey}
Yao, L., Chu, Z., Li, S., Li, Y., Gao, J., and Zhang, A. (2021).
\newblock A survey on causal inference.
\newblock {\em ACM Transactions on Knowledge Discovery from Data}, 15(5):1--46.

\bibitem[Zhang and Hyv{\"a}rinen, 2009]{zhang2009postnonlinear}
Zhang, K. and Hyv{\"a}rinen, A. (2009).
\newblock On the identifiability of the post-nonlinear causal model.
\newblock In {\em Proceedings of the 25th Conference on Uncertainty in Artificial Intelligence (UAI)}, pages 647--655. AUAI Press.

\bibitem[Zhou et~al., 2024]{zhou2024promises}
Zhou, Y., Tang, D., Kong, D., and Wang, L. (2024).
\newblock Promises of parallel outcomes.
\newblock {\em Biometrika}, 111(2):537--550.

\end{thebibliography}
% JDS often provides its own bst; if required, switch to jds style:
% \bibliographystyle{jds}
\bibliographystyle{apalike}

\end{document}